\documentclass[a4paper]{article}
\usepackage[english]{babel}
\usepackage[utf8x]{inputenc}
\usepackage[T1]{fontenc}
\newlength\figureheight
\newlength\figurewidth

\usepackage[a4paper,top=3cm,bottom=2cm,left=2cm,right=2cm,marginparwidth=1.75cm]{geometry}

\usepackage{amsmath}
\usepackage{graphicx}
\usepackage[colorlinks=true, allcolors=blue]{hyperref}
\usepackage{xcolor}
\usepackage{amsthm}
\usepackage{amssymb}
\usepackage{enumitem}   
\usepackage[toc,page]{appendix}

\numberwithin{equation}{section}

\DeclareMathOperator*{\argmin}{argmin}

\newtheorem{cor}{Corollary}[section]
\newtheorem{thm}{Theorem}[section]
\newtheorem{definition}{Definition}[section]
\newtheorem{lemma}{Lemma}[section]

\newcommand{\R}{\mathbb R}
\newcommand{\Hu}{\mathbf H_1}
\newcommand{\Hd}{\mathbf H_2}
\newcommand{\smallO}[1]{\ensuremath{\mathop{}\mathopen{}{\scriptstyle\mathcal{O}}\mathopen{}\left(#1\right)}}
\newcommand{\bigO}[1]{\ensuremath{\mathop{}\mathopen{}{\large\mathcal{O}}\mathopen{}\left(#1\right)}}

\newcommand{\dpy}{\displaystyle}

\title{\textbf{Nesterov's acceleration and Polyak's heavy ball method in continuous time: convergence rate analysis under geometric conditions and perturbations}}

\author{Othmane Sebbouh$^1$, Charles Dossal$^1$, Aude Rondepierre$^{1,2}$\\
\mbox{} \\
$^1$ IMT, Univ. Toulouse, INSA Toulouse, France. \\
$^3$ LAAS, Univ. Toulouse, CNRS, Toulouse, France.
\mbox{} \\
{\scriptsize othmane.sebbouh@gmail.com,
$\{$Charles.Dossal,Aude.Rondepierre$\}$@insa-toulouse.fr}
}
\date{}
\begin{document}
\maketitle
\begin{abstract}
 In this article a family of second order ODEs associated to inertial gradient descend is studied. 
These ODEs are widely used to build trajectories converging to a minimizer $x^*$ of a function $F$, possibly convex.
This family includes the continuous version of the Nesterov inertial scheme and the continuous heavy ball method.  
Several damping parameters, not necessarily vanishing, and a perturbation term $g$ are thus considered. 
The damping parameter is linked to the inertia of the associated inertial scheme and the perturbation term $g$ 
is linked to the error that can be done on the gradient of the function $F$. 
This article presents new asymptotic bounds on $F(x(t))-F(x^*)$ where $x$ is a solution of the ODE, when $F$ is convex and satisfies
local geometrical properties such as \L ojasiewicz properties and under integrability conditions on $g$. 
Even if geometrical properties and perturbations were already studied for most ODEs of these families, it is the first time they are jointly studied. 
All these results give an insight on the behavior of these inertial and perturbed algorithms if $F$ satisfies some \L ojasiewicz properties especially in the setting 
of stochastic algorithms.
\end{abstract}


\tableofcontents

\section{Introduction}

Let $F:\R^n\rightarrow \R$ be a differentiable convex function admitting at least one minimizer. In this paper we study the asymptotic behavior of the trajectories of the perturbed second-order ordinary differential equation (ODE):
\begin{equation} \label{ODE}
\ddot{x}(t) + \beta(t)\dot{x}(t) + \nabla	F(x(t)) = g(t),
\end{equation}
where $t_0>0$, $\beta(t)= \dfrac{\alpha}{t^\theta}$, with $\alpha>0$ and $\theta\in [0,1]$, is a viscous damping coefficient and $g:[t_0,+\infty[\rightarrow \R$ an integrable source term that can be interpreted as a small external perturbation exerted on the system. Throughout the paper, we assume that, for any initial conditions ($x_0,v_0)\in \R^n\times\R^n$, the Cauchy problem associated with the differential equation \eqref{ODE}, has a unique global solution satisfying $(x(t_0),\dot x(t_0) )=(x_0,v_0)$. This is guaranteed for instance when the gradient function $\nabla F$ is Lipschitz  on bounded subsets of $\R^n$ \cite{HarauxJendoubi2012,jendoubi2015asymptotics}.

During the last five years many articles study these ODEs, the convergence of the trajectory $x(t)$ or the decay rate of $F(x(t))$ to its minimum value $F^*$, see for example \cite{attouchpeypouquet2016,su2016differential,may2015asymptotic} and reference therein. In \cite{su2016differential} Su et al. proved the Nesterov acceleration scheme can be seen as a discretization scheme of the ODE \eqref{ODE} with $\theta=1$ and $g(t)=0$. Moreover the convergence properties of the solution $x(t)$ of \eqref{ODE} are directly linked with the ones of the sequence defined by the Nesterov scheme and the Lyapunov analysis used in both cases to prove the convergence are very similar. As another example, the choice $\theta=0$ corresponds to the heavy ball damping. It turns out that this family of ODEs is related to inertial optimization algorithms, with various inertia, depending on the choice of the damping function $\frac{\alpha}{t^\theta}$ and including perturbation or error terms defined by $g$. Many results concerning inertial algorithms have been transposed to the continuous setting such as the convergence of FISTA iterates in \cite{chambolle2015convergence} by Chambolle et al. which has been transposed by Attouch et al. \cite{attouch2018fast} to the weak convergence of the trajectory of the solution of \eqref{ODE} with $\theta=1$ and $\alpha>3$. Conversely 
May in \cite{may2015asymptotic} and Attouch et al. in \cite{attouch2018fast} proved that for $\theta=1$ and $\alpha>3$, if $F$ is convex, the solution $x$ of \eqref{ODE} satisfies $F(x(t))-F^*=o\left(\frac{1}{t^2}\right)$ and this result has been extended to the sequence generated by FISTA by Attouch et al. in \cite{attouchpeypouquet2016}.
Consequently, studying \eqref{ODE} is also a first step to have a better understanding of general and perturbed inertial schemes to minimize convex functions.  

In \cite{cabot2009long} Cabot et al. consider a general damping term and a vanishing perturbation term $g=0$. Their study gives decay rates on $F(x(t))-F^*$ when $\theta\in [0,1)$. If only a convexity assumption is made on $F$, Su et al \cite{su2016differential} proved  that if $\theta=1$ and $\alpha\geqslant 3$, we can get 
$F(x(t))-F^*=O\left(\frac{1}{t^2}\right)$. In \cite{attouch2018fast,attouch2017asymptotic} authors complete these first results for $\theta=1$, when $\alpha>3$, proving the weak convergence of trajectory $x$ and showing that  $F(x(t))-F^*=o\left(\frac{1}{t^2}\right)$. In \cite{aujoldossal2017,attouch2017rate} authors give some optimal bound on $F(x(t))-F^*$ in the subcritical case $\alpha<3$. More general damping functions $\beta(t)$ have been studied by Cabot et al., Jendoubi et al., Attouch et al. see \cite{cabot2009second,jendoubi2015asymptotics,attouch2017asymptotic} for complete results. In particular, if $F$ is convex and $\theta\in[0,1)$, we can get: $F(x(t))-F^*=O\left(\frac{1}{t^{1+\theta}}\right)$.

Several works extend these previous results with a non vanishing perturbation term $g$ proposing some integrability conditions on $g$, see for example Balti et al. \cite{balti2016asymptotic} for $\theta\in[0,1)$ and Attouch et al. 
\cite{attouch2018fast} for $\theta=1$ and $\alpha\geqslant 3$. In these two settings, the condition on $g$ ensuring the optimal decay rate 
$F(x(t))-F^*=O\left(\frac{1}{t^{1+\theta}}\right)$ is the following:
\begin{equation}\label{IntCond1}
\int_{t_0}^{+\infty}t^{\frac{1+\theta}{2}}\|g(t)\|dt<+\infty.
\end{equation}
For $\theta=1$ and $\alpha<3$, Attouch et al. and Aujol et al. \cite{attouch2017rate,aujoldossal2017}  proved that this condition can be weakened to 
\begin{equation}
\int_{t_0}^{+\infty}t^{\frac{\alpha}{3}}\|g(t)\|dt<+\infty
\end{equation}
to ensure that $F(x(t))-F^*=O\left(\frac{1}{t^{\frac{2\alpha}{3}}}\right)$.
%
%
In \cite{attouch2017asymptotic} Attouch et al. for $\theta\in(0,1]$ and in \cite{aujol2018optimal} Aujol et al. for $\theta=1$ proved that these decay rates can be improved if more geometrical properties are known on $F$ when the perturbation term $g$ vanishes. These geometrical properties describe the growth of $F$ around the set of minimizers and are linked with \L ojasiewicz properties when $F$ is convex. 

The goal of this work is to generalize all the previous works providing accurate rates on $F(x(t))-F^*$ for any $\alpha>0$, for any $\theta\in[0,1]$, depending on the geometrical properties of $F$ such as \L ojasiewicz properties and integrability conditions on $g$. To our best knowledge, this is the first work combining geometrical properties on $F$ and integrability on $g$ to provide decays on $F(x(t))-F^*$. More precisely, we will always consider that $F$ is convex, has a unique minimizer and we always assume integrability conditions on $g$ that ensure the convergence of $F(x(t))-F^*$ to $0$. Consequently, the convergence of the trajectory $(x(t))_{t\geqslant t_0}$ to the unique minimizer is always ensured. That is why, in all theorems, the geometrical assumptions are only made on a neighborhood of the minimizer and are not necessarily global.

The paper is organized as follows.
In Section~\ref{sec_geom}, we introduce the geometrical hypotheses we consider on the function $F$, and their relation with \L ojasiewicz property. We then present the contributions of the paper in Section~\ref{sec_contrib}: depending on the geometry of the function $F$ and the value of the damping parameters $\alpha$ and $\theta$, we show that combining a flatness condition and a sharpness condition such as the \L ojasiewicz property provides new and better convergence rates for the values $F(x(t))-F^*$. The proofs of the theorems are given in Section~\ref{sec_proofs}. Some technical proofs are postponed to Appendix~\ref{sec_appendix}.

\section{Preliminaries: local geometry of convex functions}\label{sec_geom}

In this section we recall some definitions and results concerning the local geometry of convex functions around their set of minimizers, see \cite{aujol2018optimal} for more details. 

Throughout the paper, we assume that the ODE \eqref{ODE} is defined in $\R^n$ equipped with the euclidean scalar product $\langle\cdot,\cdot\rangle$ and the associated norm $\|\cdot\|$. As usual $B(x^*,r)$ denotes the open euclidean ball with center $x^*\in \R^n$ and radius $r>0$. We now introduce on the one hand a flatness assumption that ensures that the function is not too sharp in the neighborhood of its minimizers, and on the other hand a sharpness assumption ensuring that the magnitude of the gradient is not too low in the neighborhood of the minimizers. 

\begin{definition}
Let $F : \mathbb{R}^n \rightarrow \mathbb{R}$ be a convex differentiable function with $X^{*} = \argmin F \neq \emptyset$, and $F^{*} = \inf f$. 
\begin{enumerate}
\item Let $\gamma \geqslant 1$. The function $F$ satisfies the condition $\mathbf{H_{1}}(\gamma)$ if, for any minimizer $x^* \in X^* $, there exists $\eta > 0$ such that:
\begin{equation*}
\forall x \in B(x^*, \eta),~F(x) - F^* \leqslant \\  \frac{1}{\gamma} \langle \nabla F(x), x - x^* \rangle.
\end{equation*}
\item Let $r \geqslant 1$. The function $F$ satisfies the growth condition $\mathbf{H_{2}}(r)$ if for any minimizer $x^* \in X^* $, there exist $K_r > 0$ and $\epsilon > 0$ such that:
\begin{equation*}
\forall x \in B(x^*, \epsilon), \  K_r \text{d}(x, X^*)^r \leqslant F(x) - F^*.
\end{equation*}
\end{enumerate}
\end{definition}
The assumption $\Hu(\gamma)$ has been already used in \cite{cabot2009long,su2016differential,aujoldossal2017,aujol2018optimal,apidopoulos2018discrete}. Note that any convex differentiable function satisfies $\Hu(1)$ and that any differentiable function such that $(F - F^*)^{\frac{1}{\gamma}}$ is convex for some $\gamma\geqslant 1$, satisfies $\mathbf{H_{1}}(\gamma)$. More precisely, the hypothesis $\Hu(\gamma)$ can be seen as a flatness condition on the geometry of a convex function around its sets of minimizers \cite[Lemma 2.4]{aujol2018optimal}: any convex differentiable function $F$ satisfying $\Hu(\gamma)$ for some $\gamma\geqslant 1$, also satisfies: for any minimizer $x^*\in X^*$, there exist $M>0$ and $\eta>0$ such that:
\begin{equation}
\forall x\in B(x^*,\eta),~F(x) -F(x^*) \leq M \|x-x^*\|^\gamma.\label{cd:flatness}
\end{equation}

The hypothesis $\Hd(r)$ with $r\geq 1$, is a growth condition on the function $F$ around its set of minimizers (critical points in the non-convex case) ensuring that $F$ is sufficiently sharp (at least as sharp as  $x\mapsto \|x-x^*\|^r$) in the neighborhood of $X^*$. It is also called r-conditioning \cite{garrigos2017convergence} or H\"olderian error bounds \cite{Bolte2017}. In the convex setting, this growth condition is equivalent to the \L ojasiewicz inequality \cite{Loja63,Loja93}, a key tool in the mathematical analysis of continuous and discrete dynamical systems, with exponent $\theta = 1-\frac{1}{r}\in (0,1]$:
\begin{definition}
\label{def_loja}
A differentiable function $F:\R^n \to \mathbb R$ is said to have the \L ojasiewicz property with exponent $\theta \in [0,1)$ if, for any critical point $x^*$, there exist $c> 0$ and $\varepsilon >0$ such that:
\begin{equation*}
\forall x\in B(x^*,\varepsilon),~\|\nabla F(x)\| \geqslant c\left(F(x)-F(x^*)\right)^{\theta}.\label{loja}
\end{equation*}
where: $0^0=0$ when $\theta=0$ by convention.
\end{definition}

Typical examples of functions having the \L ojasiewicz property are real-analytic functions and $C^1$ subanalytic functions, or semi-algebraic functions \cite{Loja63,Loja93}. Strongly convex functions satisfy a global \L ojasiewicz property with exponent $\theta=\frac{1}{2}$ \cite{attouch2009convergence}, or equivalently a global version of the growth condition, namely:
$$\forall x\in \R^n,  F(x)-F^*\geqslant \frac{\mu}{2}d(x,X^*)^2,$$
where $\mu>0$ denotes the parameter of strong convexity. Likewise, convex functions having a strong minimizer in the sense of \cite[Section 3.3]{attouch2018convergence}, also satisfy a global version of $\Hd(2)$. By extension, uniformly convex functions of order $p\geqslant 2$ satisfy the global version of the hypothesis $\Hd(p)$ \cite{garrigos2017convergence}.

Finally, observe that any convex differentiable function $F$ satisfying both hypothesis $\Hu(\gamma)$ and $\Hd(r)$, has to be at least as flat as $\|x-x^*\|^\gamma$ and as sharp as $\|x-x^*\|^r$ in the neighborhood of its minimizers. More precisely, combining  \eqref{cd:flatness} and $\Hd(r)$, we have:
\begin{lemma}[{\cite[Lemma 2.5]{aujol2018optimal}}]
If a convex differentiable function $F$ satisfies both $\Hu(\gamma)$ and $\Hd(r)$,  with $\gamma,r\geqslant 1$, then necessarily: $r\geqslant \gamma$.\label{lem:geometry}
\end{lemma}

\section{Contributions}\label{sec_contrib}

In this section, we state convergence rates for the values $F(x(t))-F^*$ along the trajectory $x(t)$ solution of \eqref{ODE}, depending on geometrical properties $\Hu$ and $\Hd$ of the function $F$ and on integrability conditions on the perturbation term $g$. Geometry and perturbations have been studied separately in several papers for this family of ODEs, and the specificity of this work is to study both aspects jointly.

Consider first the case when $\theta =1$. If $F$ only satisfies $\Hu(\gamma)$ for some $\gamma\geqslant 1$ then for low friction parameter $\alpha$, we have the following result:
\begin{thm}[{\cite[Theorem 2]{aujoldossal2017}}] \label{th:Nesterov:sharp}
Let $\alpha > 0$ and $t_0>0$. Let $x(\cdot)$ be any solution of the ODE \eqref{ODE} with $\theta=1$ and $(x(t_0),\dot x(t_0))=(x_0,v_0)$. Assume that: 
$
\int_{t_0}^{+\infty} t^{\frac{\gamma \alpha}{\gamma + 2}}\|g(t)\|dt < +\infty.
$
If $F$ satisfies $\Hu(\gamma)$ for some $\gamma \geqslant 1$ and if $\alpha\leqslant 1+\frac{2}{\gamma}$ then:
\begin{equation}
F(x(t)) - F^* = \bigO{t^{-\frac{2\gamma \alpha}{\gamma + 2}}}.
\end{equation}
\end{thm}
This result has been first stated and proved in the unpublished report \cite[Theorem 2]{aujoldossal2017} by Aujol and Dossal in 2017 for convex differentiable functions satisfying $(F-F^*)^{\frac{1}{\gamma}}$ convex. 

For large friction parameters $\alpha$, the sole assumption $\Hu(\gamma)$ on $F$ is not sufficient anymore to obtain a decay faster than $\bigO{\frac{1}{t^2}}$ which is the uniform rate that can be achieved for $\alpha \geqslant 3$ \cite{su2016differential}. The contribution of this paper is to show that a flatness condition $\Hu$ associated to a sharpness condition such as the \L ojasiewicz property provides new and better convergence rates for the values $F(x(t))-F^*$. We can thus compare these results with classical bounds that can be achieved with geometrical assumptions on $F$, such as convexity or \L ojasiewicz properties without any perturbation term or with results dealing with a non vanishing perturbation term $g$ but with simple assumptions on $F$.


\subsection{Convergence rates for sharp geometries}
In this section, we state convergence rates on the values $F(x(t))-F^*$ along the trajectory $x(t)$, that can be achieved for functions satisfying geometrical hypothesis such as $\Hu(\gamma)$ and/or $\Hd(2)$. The cases $\theta=1$ and $\theta\in [0,1)$ are treated separately.

Let us first consider the case when $\theta=1$ i.e. the ODE:
\begin{equation} \label{ODE:Nesterov}
\ddot{x}(t) + \frac{\alpha}{t}\dot{x}(t) + \nabla	F(x(t)) = g(t).
\end{equation} 
\begin{thm} \label{th:Nesterov:sharp2}
Let $\alpha > 0$ and $t_0>0$. Let $x(\cdot)$ be any solution of the ODE \eqref{ODE:Nesterov} with $\theta=1$ and $(x(t_0),\dot x(t_0))=(x_0,v_0)$. Assume that: 
\begin{equation}
\int_{t_0}^{+\infty} t^{\frac{\gamma \alpha}{\gamma + 2}}\|g(t)\|dt < +\infty.
\end{equation}
If $F$ satisfies $\Hu(\gamma)$ and $\Hd(2)$, for some $\gamma \leqslant 2$, if $F$ has a unique minimizer and if $\alpha > 1 + \frac{2}{\gamma}$, then
\begin{equation}
F(x(t)) - F^* = \bigO{t^{-\frac{2\gamma \alpha}{\gamma + 2}}}.
\end{equation}
\end{thm}

Note that Theorem~\ref{th:Nesterov:sharp2} whose proof is detailed in Section \ref{proof:Nesterov:sharp}, only applies for $\gamma\leqslant 2$ since according to Lemma \ref{lem:geometry}, there exists no function satisfying both $\Hu(\gamma)$ and $\Hd(2)$ for $\gamma>2$. Moreover the integrability condition given in Theorem~\ref{th:Nesterov:sharp2} generalizes the integrability condition given in \cite{aujoldossal2017} to any $\alpha >0$ under the growth condition $\Hd(2)$ and coincides in the limit case $\alpha=3$, which was expected. Theorems \ref{th:Nesterov:sharp2} can be seen as an extension of former results with a non vanishing perturbation term $g$, see \cite[Theorem 4.2]{aujol2018optimal}. 

Let $\alpha > 0$ and $\theta\in [0,1)$. We now consider the heavy ball system with a general friction term:
\begin{equation} \label{ODE:HB}
\ddot{x}(t) + \frac{\alpha}{t^\theta}\dot{x}(t) + \nabla	F(x(t)) = g(t).
\end{equation} 
\begin{thm} \label{th:HB:sharp}
Let $\gamma \in [1,2]$, $m\in \left(0,\frac{2\gamma}{\gamma+2}\right)$ and $t_0>0$. Note: $\Gamma(t) = \int_{t_0}^{t} \frac{\alpha}{s^\theta}ds$. Let $x(\cdot)$ be any solution of the ODE \eqref{ODE:HB} with $\theta\in [0,1)$ and $(x(t_0),\dot x(t_0))=(x_0,v_0)$. Assume that: $$\int_{t_0}^{+\infty} e^{m\Gamma(t)}\|g(t)\|dt < +\infty.$$ 
If $F$ satisfies $\Hu(\gamma)$ and $\Hd(2)$, and admits a unique minimizer $x^*$ then:
\begin{eqnarray*}
F(x(t)) - F^* = \bigO{e^{-m\Gamma(t)}},\quad \|x(t)-x^*\|^2 = \bigO{e^{-m\Gamma(t)}},\quad\|\dot x(t)\|^2= \bigO{e^{-m\Gamma(t)}}.
\end{eqnarray*}
\end{thm}
Theorem~\ref{th:HB:sharp} when $\theta\in(0,1)$ can be seen as an extension of \cite[Theorem 6.1]{attouch2017asymptotic} to functions with some geometrical properties as $\Hu$ and $\Hd$. Its proof is detailed in Section~\ref{sec_proofs} and is an extension of the proof of \cite[Theorem 3.12]{attouch2017asymptotic} to a non-vanishing perturbation term $g\neq 0$. Note that \cite{attouch2017asymptotic} deals only with vanishing damping, that is $\theta>0$, while Theorem \ref{th:HB:sharp} deals also with the case $\theta=0$. 

Observe also that Theorem \ref{th:HB:sharp} also applies in the case $\theta=1$ and provides convergence rates in $\bigO{t^{-m\alpha}}$ for any $m\in \left(0,\frac{2\gamma}{\gamma+2}\right)$, which is slower but infinitely close to the convergence rate $\bigO{t^{-\frac{2\gamma \alpha}{\gamma + 2}}}$ provided by Theorem~\ref{th:Nesterov:sharp2}.

Finally, observe that the integrability conditions given in these three theorems are always stronger than the condition $\int_{t_0}^{+\infty} t^p\|g(t)\|dt<+\infty$, with $p=\min(1,\frac{\alpha}{3})$, given in \cite[Theorem 5.1]{attouch2017rate} and ensuring the convergence of the values $F(x(t))-F^*$ to $0$. Consequently the trajectory $x(t)$ actually converges to the unique minimizer of $F$ so that the geometrical assumptions $\Hu$ and $\Hd$ can be used locally that is in the neighborhood of the unique minimizer of $F$.

\subsection{Convergence rates for flat geometries}
In this section, we state new convergence rates on the values $F(x(t))-F^*$ along the trajectory $x(t)$, that can be achieved for functions satisfying geometrical hypothesis such as $\Hu(\gamma)$ and $\Hd(\gamma)$ for any $\gamma >2$. The cases $\theta=1$ and $\theta\in [0,1)$ are treated jointly.

Let us first consider the unperturbed case ($g=0$). We announce new results on the convergence of the function values along the trajectory $x(t)$ with additional geometrical assumptions, but without perturbations:
\begin{thm}
\label{th:unperturbed:flat case}
Let $\gamma_1 > 2$, $ \gamma_2\geqslant \gamma_1 $. Note $r = \frac{1 + \theta}{2}$. Suppose $x$ is a solution to the ODE (\ref{ODE}) with $g=0$. If $F$ is coercive and satisfies $\mathbf{H_1}(\gamma_1)$ and $\mathbf{H_2}(\gamma_2)$, and 
\begin{enumerate}
    \item if $\theta=1$ and $\alpha\geqslant \frac{\gamma_1+2}{\gamma_1 - 2}$ or 
    \item if $\theta<1$,
\end{enumerate}
then:
\begin{equation}
F(x(t)) - F^* = \bigO{\frac{1}{t^{\frac{2r\gamma_2}{\gamma_2 - 2}}}}.
\end{equation}
\end{thm}
Note that Theorem~\ref{th:unperturbed:flat case} is a generalization of Theorem 4.3 in \cite{aujol2018optimal} available for $\theta=1$ to any $\theta\in(0,1]$. We now prove that the convergence rates provided by Theorem~\ref{th:unperturbed:flat case} remain valid in the perturbed case ($g\neq 0$):
\begin{thm}\label{th:perturbed:flat case}
Let $\alpha >0$, $\theta\in [0,1]$ and $t_0>0$. Let $x$ be the solution of the ODE \eqref{ODE} for given initial conditions $(x(t_0),\dot x(t_0))=(x_0,v_0)$. Let $\gamma_1 > 2$, $ \gamma_2\geqslant \gamma_1$ and $r = \frac{1 + \theta}{2}$. 
Assume that:
\begin{equation}
\int_{t_0}^{+\infty} t^{\frac{r \gamma_2}{\gamma_2 - 2}}\|g(t)\|dt < +\infty.
\end{equation}
If $F$ satisfies $\mathbf{H_1}(\gamma_1)$ and $\mathbf{H_2}(\gamma_2)$ and admits a unique minimizer then:
\begin{enumerate}
    \item if $\theta=1$ and $\alpha\geqslant \frac{\gamma_1+2}{\gamma_1 - 2}$ or 
    \item if $\theta<1$,
\end{enumerate}
then we have
\begin{equation}
F(x(t)) - F^* = \bigO{\frac{1}{t^{\frac{2r \gamma_2}{\gamma_2 - 2}}}}.
\end{equation}
\end{thm}
As in \cite{aujol2018optimal} with a non vanishing perturbation term $g$, if $\gamma_1 = \gamma_2$, we have furthermore the convergence of the trajectory:
\begin{cor}
Let $\alpha >0$, $\theta\in [0,1]$ and $t_0>0$. Let $x$ be the solution of the ODE \eqref{ODE} for given initial conditions $(x(t_0),\dot x(t_0))=(x_0,v_0)$.
Let $\gamma > 2$. Note $r=\frac{1+\theta}{2}$.
Assume that:
$$\int_{t_0}^{+\infty} t^{\frac{r\gamma}{\gamma - 2}}\|g(t)\|dt < +\infty.$$
If $F$ satisfies $\mathbf{H}_1(\gamma)$ and $\mathbf{H}_2(\gamma)$ and admits a unique minimizer then:
\begin{enumerate}
    \item if $\theta=1$ and $\alpha\geqslant \frac{\gamma+2}{\gamma - 2}$ or 
    \item if $\theta<1$,
\end{enumerate}
then we have:
\begin{equation}
\|\dot{x}(t)\| = \bigO{\frac{1}{t^{\frac{r\gamma}{\gamma - 2}}}}.
\end{equation}\label{cor:flat}
\end{cor}

Note that in the case of the classical heavy ball ($\theta=0$), Theorem~\ref{th:perturbed:flat case} can be seen as an extension of \cite[Corollary 5.1]{Begout2015} using a different approach: indeed in \cite{Begout2015}, the authors proved a similar convergence rate under \L ojasiewicz properties, but without any convexity assumption on $F$. 

Observe also that to deal with a non-vanishing perturbation term, the uniqueness of the minimizer seems to be crucial despite the fact we can avoid this assumption when $g=0$, see \cite[Theorem 4.5 and Corollary 4.6]{aujol2018optimal}. 

Finally observe that in the integrability condition given in Theorem \ref{th:perturbed:flat case}, the exponent is $\frac{r\gamma_2}{\gamma_2-2}=\frac{(1+\theta)\gamma_2}{2(\gamma_2-2)}$.
Since $\frac{r\gamma_2}{\gamma_2-2} \geqslant r$, the integrability condition \eqref{IntCond1} from \cite{may2015asymptotic} is automatically satisfied and ensures that the trajectory $x(t)$ converges to the unique minimizer of $F$. The geometrical assumptions $\Hu$ and $\Hd$ thus can be used locally. Moreover:
$$\underset{\gamma_2\in(2+\infty)}{\inf}\frac{(1+\theta)\gamma_2}{2(\gamma_2-2)}=\frac{1+\theta}{2},$$ 
so that we get the same exponent as in the integrability condition in \eqref{IntCond1} which is also what we expected.


\subsection{Strategies of proofs}
The guideline of the proofs of our results is the same: Lyapunov functions (or energies). In this section we present the state of the art strategies using Lyapunov functions and a sketch of the strategies used in this paper.

Note that in the four Theorems \ref{th:Nesterov:sharp}, \ref{th:Nesterov:sharp2}, \ref{th:HB:sharp} and \ref{th:perturbed:flat case} dealing with a non vanishing perturbation term, we modified Lyapunov functions used when $g=0$, and proposed integrability conditions on $g$ that ensure the same decay that the one that can be achieved with $g=0$. That is the reason why the Lyapunov functions we use are closed to those that can be found in the literature. Note also that the Gr\"onwall-Bellman Lemma is a key lemma in each proof, but the exact way to deal with a non vanishing perturbation term $g$ is different in each theorem.  

\subsubsection{State of the art strategies}
To prove each decay of $F(x(t))-F^*$, the main idea is to define a function of $t$ which will be denoted by $\mathcal{E}$, $\mathcal{H}$ or $\mathcal{G}$, which involves the term $F(x(t))-F^*$ and which is bounded. The choice of the Lyapunov function depends on the ODE 
($\alpha$ and $\theta$) and on assumptions on function $F$ (the flatness hypothesis $\mathbf{H}_1$). A simple Lyapunov function to study solutions of \eqref{ODE} when $F$ is convex with $g=0$ is: 
\begin{equation}
\mathcal{E}(t)=F(x(t))-F^*+\frac{1}{2}\|\dot{x}(t)\|^2.
\end{equation}
Indeed, $\mathcal{E}$ is a sum of positive terms and $\mathcal{E}'(t)=-\beta(t)\|\dot{x}(t)\|^2\leqslant 0$. This simple Lyapunov function ensures that $\mathcal{E}$ is non increasing which implies that 
$F(x(t))$ is bounded. Many Lyapunov have been proposed to study \eqref{ODE:Nesterov} when $F$ is convex, see for example \cite{su2016differential,attouch2016fast,attouch2017rate,aujoldossal2017} or for more general friction term, see \cite{cabot2009long,balti2016asymptotic}. A simple example to study the specific case of Nesterov damping (i.e. $\theta=1$) when $F$ is convex and $g=0$,  is:
\begin{equation}
\mathcal{E}(t)=t^2(F(x(t))-F^*)+\frac{1}{2}\|(\alpha-1)(x(t)-x^*)+t\dot x(t)\|^2_2.
\end{equation}
Indeed, a simple calculation shows that:
\begin{equation}
\mathcal{E}'(t)\leqslant -\alpha t\|\dot{x}(t)\|^2+(3-\alpha)t(F(x(t))-F^*).
\end{equation}
From this Lyapunov function we deduce that if $\alpha\geqslant 3$, $\mathcal{E}$ is non increasing and thus that 
\begin{equation}
F(x(t))-F^*\leqslant \frac{\mathcal{E}(t_0)}{t^2}.
\end{equation}
In \cite{aujol2018optimal}, we propose to extend this Lyapunov approach to deal with geometrical properties of $F$ for the Nesterov damping i.e \eqref{ODE:Nesterov}  with $g=0$ using Lyapunov functions 
$\mathcal{H}$ :
\begin{equation}
\mathcal{H}(t)=t^p(F(x(t))-F^*)+R(x(t)),
\end{equation}
where $R$ is non necessarily positive and $p$ depends on properties of the damping parameter $\alpha$ and on the geometric properties of $F$ to get $F(x(t))-F^*=O\left(\frac{1}{t^p}\right)$.
To get such a bound we first bound the energy $\mathcal{H}$ and thus use the \L ojasiewicz properties of $F$ ($\mathbf{H}_2$) to deduce the bound on $F(x(t))-F^*$.

For other choices of function $\beta(t)$, especially when $\theta<1$ in \eqref{ODE}, the decay may be faster than polynomial if $F$ satisfies some \L ojasiewicz properties, see \cite{polyak2017lyapunov,attouch2017asymptotic,Bolte2017}. A way to prove these faster decay is to build Lyapunov energies satisfying some differential inequalities for suitable function $\gamma$ :
\begin{equation}\label{Ineqdiff}
\mathcal{E}'(t)\leqslant -\gamma(t)\mathcal{E}(t)
\end{equation}  
which implies that 
\begin{equation}
\mathcal{E}(t)\leqslant \mathcal{E}(t_0)e^{-\int_{t_0}^t\gamma(s)ds}.
\end{equation}
Moreover, a simple way to deal with the perturbation term $g$ in \eqref{ODE} is to add an integral term in the Lyapunov energy $\mathcal{E}$ or $\mathcal{H}$ depending on $x$ and $g$, see for example \cite{balti2016asymptotic,attouch2016fast,aujoldossal2017} and references therein and to use a Gr\"onwall-Bellman lemma to conclude.

\subsubsection{Sketch of proof}
For each theorem, we define a Lyapunov Energy  $\mathcal{G}$ (or $\mathcal{H}$) defined by parameters that are set depending on the hypotheses of each theorem. This function $\mathcal{G}$ also depends on the perturbation term $g$. The first step of the proof consists in proving that $\mathcal{G}$ is bounded. The second step use a Gr\"onwall-Bellman lemma  and hypothesis $\mathbf{H}_2$ to conclude. 
More precisely, Theorems \ref{th:Nesterov:sharp2} and \ref{th:HB:sharp} deal with sharp functions, that is functions satisfying $\mathbf{H}_2(2)$, or \L ojasiewicz properties with parameter equal to $\frac{1}{2}$.  The first one is dedicated to Nesterov damping ($\theta=1$). The polynomial rate was known in the case where $g=0$, see \cite{aujol2018optimal}, and we propose to modify the Lyapunov function that is used by adding an integral term :
\begin{equation*}
\mathcal{G}(t)=t^p\left(t^2(F(x(t))-F^*)+\frac{1}{2}\|\lambda(x(t)-x^*)+t\dot{x}(t)\|_2^2+\frac{\xi}{2}\|x(t)-x^*\|^2\right)+\int_{t}^Ts^{\delta}\langle \lambda (x(s)-x^*)+s\dot{x}(s),g(s)\rangle ds
\end{equation*}
and prove that $\mathcal{G}$ is bounded using differential inequalities. The Gr\"onwall-Bellman lemma and integrability hypotheses on $g$ are used to conclude. 

Theorem  \ref{th:HB:sharp} deals with the case $\theta<1$. In this case it is known \cite{polyak2017lyapunov,attouch2017asymptotic,Bolte2017} that the decay of $F(x(t))-F^*$
 is faster than polynomial. We propose to use a Lyapunov function similar to the previous one. This time $\lambda$
 and $\xi$ may be functions of $t$ : 
 \begin{equation*}
\mathcal{G}(t)=F(x(t))-F^*+\frac{1}{2}\|\lambda(t)(x(t)-x^*)+\dot{x}(t)\|_2^2+\frac{\xi(t)}{2}\|x(t)-x^*\|^2+\int_{t}^Ts^{\delta}\langle \lambda (x(s)-x^*)+s\dot{x}(s),g(s)\rangle ds.
\end{equation*}
For suitable choices of $\lambda, \xi$ and $\delta$ we can get some differential inequalities like \eqref{Ineqdiff} satisfied by $\mathcal{G}$  and conclude using Gr\"onwall-Bellman lemma.

Theorems \ref{th:unperturbed:flat case} and \ref{th:perturbed:flat case} deal with {\it flat functions} that is functions $F$ satisfying $\mathbf{H}_2(\gamma_2)$ with $\gamma_2>2$, i.e 
\L ojasiewicz properties with an exponent greater than $\frac{1}{2}$. Theorem  \ref{th:unperturbed:flat case} focuses on the case $g=0$ and Theorem  \ref{th:perturbed:flat case} is the general perturbed case. Both theorems provide results for $\theta\in[0,1]$ including Heavy Ball ($\theta=0$) and Nesterov ($\theta=1$). In both theorems, the bound given on 
$F(x(t))-F^*$ is polynomial. For Theorem  \ref{th:unperturbed:flat case}, inspired by \cite{aujol2018optimal} we define an energy function 
\begin{equation*}
\mathcal{H}(t)=t^p\left(t^2(F(x(t))-F^*)+\frac{1}{2}\|\lambda(x(t)-x^*)+t\dot{x}(t)\|_2^2+\frac{\xi(t)}{2}\|x(t)-x^*\|^2\right)
\end{equation*}
for a suitable choice of parameters $\lambda$ and $p$, and function $\xi$. And we prove that it exists $t_1\geqslant t_0$ such that for any $t\geqslant t_1$, $\mathcal{H}'(t)\leqslant 0$. 
The function $\xi$ may be negative  but using \L ojasiewicz properties of $F$, we can deduce bounds on $F(x(t))-F^*$. In Theorem  \ref{th:perturbed:flat case}, we consider:
\begin{equation*}
\mathcal{G}(t)=t^p\left(t^2(F(x(t))-F^*)+\frac{1}{2}\|\lambda(x(t)-x^*)+t\dot{x}(t)\|_2^2+\frac{\xi(t)}{2}\|x(t)-x^*\|^2\right)+\int_{t}^Ts^{\delta}\langle \lambda (x(s)-x^*)+s\dot{x}(s),g(s)\rangle ds
\end{equation*}
and prove that $\mathcal{G}$ is bounded. Combining approaches developed in Theorem \ref{th:unperturbed:flat case}, Gr\"onwall-Bellman lemma and integrability hypotheses on $g$, we are able to conclude.

\section{Proofs}\label{sec_proofs}
In this section, we detail the proofs of Theorem~\ref{th:Nesterov:sharp2}, Theorem~\ref{th:HB:sharp}, Theorem~\ref{th:unperturbed:flat case}, Theorem~\ref{th:perturbed:flat case} and Corollary \ref{cor:flat}.
\subsection{Proof of Theorem~\ref{th:Nesterov:sharp2}}\label{proof:Nesterov:sharp}
In this section, we prove Theorem \ref{th:Nesterov:sharp2}. For a complete proof of Theorem \ref{th:Nesterov:sharp}, we refer the reader to \cite{aujoldossal2017}.

Let $x^*$ be a minimizer of $F$ and $\lambda$, $\xi$ and $T$ three real numbers. The proof of Theorem \ref{th:Nesterov:sharp2} relies on the following energy:
\begin{equation}
\mathcal{G}(t) = t^p\mathcal E(t) + \int_t^T \langle s^{\frac{p}{2}} (\lambda (x(s) - x^*) + s \dot{x}(s)), s^{\frac{p+2}{2}} g(s) \rangle ds\label{G}
\end{equation}
where the energy $\mathcal E$ is defined by:
\begin{equation}
\mathcal{E}(t) = t^2(F(x(t)) - F^*) + \frac{1}{2}\|\lambda(x(t) - x^*) + t\dot{x}(t)\|^2 + \frac{\xi}{2}\|x(t) - x^*\|^2.
\end{equation}
Using the following notations: 
\begin{equation}
\begin{aligned}\label{eqdefabc}
a(t) & = t(F(x(t)) - F^*), \\
b(t) & = \frac{1}{2t}\|\lambda(x(t) - x^*) + t\dot{x}(t)\|^2, \\
c(t) & = \frac{1}{2t}\|x(t) - x^*\|^2,
\end{aligned}
\end{equation}
we then have:
\begin{equation}
\mathcal{E}(t) = t(a(t) + b(t) + \xi c(t)).
\end{equation}

Note that the functions $\mathcal E(t)$ and $\mathcal{H}(t)=t^p\mathcal E(t)$ denote the same Lyapunov functions as those used in the proof of Theorem 4.1. in \cite{aujol2018optimal} in the non perturbed case ($g=0$). 
Our proofs are based on the following lemma:
\begin{lemma} \label{lemma Nesterov}
Let $\gamma \geqslant 1$. If $F$ satisfies the hypothesis $\mathbf{H}_1(\gamma)$ and if $\xi = \lambda(\lambda + 1 - \alpha)$, then
\begin{equation}
\mathcal{G}'(t)  \leqslant t^p\left((2 + p - \gamma \lambda)a(t) + (p + 2\lambda +2 - 2\alpha)b(t)  + \lambda(\lambda + 1 - \alpha)(p-2\lambda)c(t)\right) \label{derivation G thm 3.5.2}
\end{equation}
\end{lemma}
This lemma whose proof is detailed in Appendix \ref{Appendix Nesterov}, is the generalization of \cite[Lemma 5.1]{aujol2018optimal} to the perturbed case: the integral term in \eqref{G} was chosen to cancel the terms in $\mathcal G'(t)$ coming from the perturbation $g(t)$: following the exact same calculation steps as in the proof of \cite[Lemma 5.1]{aujol2018optimal}, we observe that all the terms coming from the perturbation $g(t)$ cancel each other out, so that we obtain the same formula as in the non perturbed version.

Choosing now $p = \frac{2\gamma \alpha}{\gamma + 2} - 2$ and $ \lambda = \frac{2 \alpha}{\gamma + 2}$ and thus:
$$\xi = \frac{2\alpha}{(\gamma+2)^2}(2+\gamma(1-\alpha)),$$ 
we get:
\begin{equation} \label{thm 1 new central}
\mathcal{G}'(t) \leqslant K_1 t^p c(t)
\end{equation}
where $K_1 = \xi(p-2\lambda) = \frac{2\xi}{\gamma+2}\left((\gamma-2)\alpha-(\gamma+2)\right)$. Since $\alpha > 1 + \frac{2}{\gamma}$ and $\gamma\leqslant 2$, we necessarily have: $\xi <0$, and thus $K_1>0$. Consequently, the energy $\mathcal E(t)$ is not a sum of non-negative terms and we cannot conclude that the energy function $\mathcal{G}$ is decreasing. To get the expected estimate on the energy, we need an additional growth condition $\Hd(2)$ to bound the term $\|x(t)-x^*\|^2$ as done in \cite{aujol2018optimal}.

Using the uniqueness of the minimizer and the fact that $F$ satisfies $\mathbf{H}_2(2)$, there exists $K > 0$ such that:
\begin{equation}
Kt\|x(t) - x^*\|^2 \leqslant t(F(x(t)) - F^*) = a(t),
\end{equation}
hence:
\begin{equation} \label{c(t):smallO}
c(t) \leqslant \frac{1}{2Kt^2}a(t).
\end{equation}
Since $\xi < 0$ with our choice of parameters, we get:
\begin{equation}
\mathcal{H}(t) \geqslant t^{p+1}(a(t) + \xi c(t)) \geqslant t^{p+1}(1 + \frac{\xi}{2Kt^2})a(t),
\end{equation}
so that there exists $t_1 \geqslant t_0$ such that for all $t \geqslant t_1$, we have: 
\begin{equation} \label{th1:ineqH}
\mathcal{H}(t) \geqslant \frac{1}{2}t^{p+1}a(t) \geqslant 0.
\end{equation}
Now from (\ref{thm 1 new central}), (\ref{c(t):smallO}) and (\ref{th1:ineqH}), we have:
\begin{equation}
\forall t\geqslant t_1,~\mathcal{G}'(t) \leqslant \frac{K_1}{K} \frac{\mathcal{H}(t)}{t^3}.
\end{equation}
Observe now that: $\mathcal{H}'(t) = \mathcal{G}'(t) + \langle  \lambda (x(t) - x^*) + t \dot{x}(t), t^{p+1} g(t) \rangle $ so that we get the following differential inequality on the energy $\mathcal H$:
\begin{equation}
\begin{aligned}
\mathcal{H}'(t) & \leqslant \frac{K_1}{K} \frac{\mathcal{H}(t)}{t^3} + \langle t^{\frac{p}{2}} (\lambda (x(t) - x^*) + t \dot{x}(t)), t^{\frac{p+2}{2}} g(t) \rangle \\
& \leqslant \frac{K_1}{K} \frac{\mathcal{H}(t)}{t^3} + t^{\frac{p}{2}}\| \lambda (x(t) - x^*) + t \dot{x}(t)\| t^{\frac{p+2}{2}}\| g(t) \| \\
&= \frac{K_1}{K} \frac{\mathcal{H}(t)}{t^3} + \sqrt[]{2}  (t^{p+1}b(t))^{\frac{1}{2}} t^{\frac{\gamma \alpha}{\gamma + 2}}\| g(t) \|
\end{aligned}
\end{equation}
Using again the fact that $c(t) =\smallO{a(t)}$ (see \eqref{c(t):smallO}), there exists $t_2\geqslant t_1$ such that, for all $t \geqslant t_2$, $a(t) + \xi c(t) \geqslant \frac{1}{2}a(t)$, which implies that:
$$\forall t\geqslant t_2,~\mathcal{H}(t)=t^{p+1}b(t) + t^{p+1}(a(t) + \xi c(t)) \geqslant \frac{1}{2}t^{p+1}a(t)+  t^{p+1}b(t)\geqslant t^{p+1}b(t).$$
Hence:
\begin{equation}
\forall t\geqslant t_2,~\mathcal{H}'(t) \leqslant \frac{K_1}{K} \frac{\mathcal{H}(t)}{t^3} + \sqrt[]{2}\mathcal{H}(t)^{\frac{1}{2}} t^{\frac{\gamma \alpha}{\gamma + 2}}\| g(t) \|.
\end{equation}
Dividing both sides of the inequality by $2\mathcal{H}(t)^{\frac{1}{2}}$ and
integrating between $t_2$ and $t$, we get:
\begin{equation}
\forall t\geqslant t_2,~\mathcal{H}(t)^{\frac{1}{2}} \leqslant \frac{K_1}{2K} \int_{t_2}^{t}  \frac{\mathcal{H}(s)^{\frac{1}{2}} }{s^3}ds + \frac{\sqrt[]{2}}{2}\int_{t_2}^{t} s^{\frac{\gamma \alpha}{\gamma + 2}}\| g(s) \|ds.
\end{equation}
Since $\int_{t_2}^{+\infty} s^{\frac{\gamma \alpha}{\gamma + 2}}\| g(s) \|ds < +\infty$:
\begin{eqnarray*}
\forall t\geqslant t_2,~\mathcal{H}(t)^{\frac{1}{2}} &\leqslant& \frac{K_1}{2K} \int_{t_2}^{t}  \frac{\mathcal{H}(s)^{\frac{1}{2}} }{s^3}ds + \frac{\sqrt[]{2}}{2}\int_{t_2}^{+\infty} s^{\frac{\gamma \alpha}{\gamma + 2}}\| g(s) \|ds\\
&\leqslant& \beta +  \int_{t_1}^{t} \frac{K_1}{2K s^3} \mathcal{H}(s)^{\frac{1}{2}}ds,
\end{eqnarray*}
where $\beta = \frac{\sqrt[]{2}}{2} \dpy\int_{t_2}^{+\infty} s^{\frac{\gamma \alpha}{\gamma + 2}}\| g(s) \|ds$. Applying the Gr\"onwall Lemma, we finally get:
\begin{equation}
\forall t\geqslant t_2,~\mathcal{H}(t) \leqslant \beta^2 \ \exp \Big(\frac{K_1}{K}\int_{t_2}^{t} \frac{1}{s^3} ds\Big) \leqslant \beta^2 \ \exp \Big(\frac{K_1}{K}\int_{t_2}^{+\infty} \frac{1}{s^3} ds\Big) = \beta^2 \exp\Big(\frac{K_1}{2Kt_2^2}\Big)
\end{equation}
In other words, we found a constant $A > 0$ such that for all $t \geqslant t_2$, $\mathcal{H}(t) \leqslant A$. According to \eqref{th1:ineqH}, we conclude that $\frac{1}{2}t^{p+2}(F(x(t)) - F^*) = \frac{1}{2}t^{p+1}a(t)$ is bounded which ends the proof of Theorem \ref{th:Nesterov:sharp2}.



\subsection{Proof of Theorem \ref{th:HB:sharp}}
Let $\lambda$, $\xi$ and $T$ three real numbers. Let $x^*$ be a minimizer of $F$ and $x(\cdot)$ any trajectory solution of:
$$\ddot{x}(t) +\beta(t)\dot{x}(t) + \nabla	F(x) = g(t),$$
where: $\beta(t)=\frac{\alpha}{t^\theta}$ with $\alpha>0$ and $\theta\in [0,1)$. The proof of Theorem \ref{th:HB:sharp} relies on the following energy:
\begin{equation}
\mathcal{G}(t) = \mathcal E(t) + \int_t^T \langle \lambda (x(s) - x^*) + \dot{x}(s)), g(s) \rangle ds\label{GHB}
\end{equation}
where:
\begin{equation}
\mathcal{E}(t) = F(x(t)) - F^* + \frac{1}{2}\|\lambda(x(t) - x^*) + \dot{x}(t)\|^2 + \frac{\xi}{2}\|x(t) - x^*\|^2.
\end{equation}

\paragraph{Case $\theta=0$.} Remember that in that case, the friction coefficient is constant: $\forall t, ~\beta(t)=\alpha$. Using the following notations: 
\begin{equation}
\begin{aligned}
a(t) & = F(x(t)) - F^*, ~b(t)  = \frac{1}{2}\|\lambda(x(t) - x^*) + \dot{x}(t)\|^2, \\
c(t) & = \frac{1}{2}\|x(t) - x^*\|^2,
\end{aligned}
\end{equation}
the energy $\mathcal E(t)$ can be rewritten as:
\begin{equation}
\mathcal{E}(t) = a(t) + b(t) + \xi c(t).
\end{equation}
The proof of Theorem \ref{th:HB:sharp} when $\theta=0$ relies on the following differential inequality whose proof in detailed in appendix:
\begin{lemma}
Let $\gamma\geqslant 1$, $\theta=0$ and $\lambda\in \R$. If $F$ satisfies the hypothesis $\Hu(\gamma)$ and $\xi=\lambda(\lambda-\alpha)$, then:
\begin{eqnarray*}
\forall t\geq t_0, ~\mathcal E'(t) &\leqslant& -\lambda\gamma a(t) +2(\lambda-\alpha)b(t) -2\lambda\xi c(t)+\langle g(t),\dot x(t)+\lambda(x(t)-x^*)\rangle\\
&\leqslant& -\lambda\gamma \left(a(t)+2\frac{\xi}{\gamma} c(t)\right) +2(\lambda-\alpha)b(t) +\langle g(t),\dot x(t)+\lambda(x(t)-x^*)\rangle.
\end{eqnarray*}\label{lem:HB:sharp}
\end{lemma}
The scheme of the rest of the proof is quite standard: we first need to control the terms in $b(t)$ and $c(t)$ in Lemma \ref{lem:HB:sharp} to deduce some differential inequality on the energy $\mathcal E$.

Let us choose $\lambda <\alpha$. In that case, the energy $\mathcal E$ is not a sum of non-negative terms anymore:
$$\mathcal{E}(t) = a(t) + b(t) + \xi c(t)$$
since: $\xi =\lambda(\lambda-\alpha) <0$. Using the growth condition $\Hd(2)$ combined by the uniqueness of the minimizer of $F$, to bound $\|x(t)-x^*\|$, we get the following inequality: there exists $t_1\geqslant t_0$ such that:
\begin{equation}
\forall t\geqslant t_1,~a(t)+2\frac{\xi}{\gamma} c(t)=a(t)-2\frac{|\xi|}{\gamma} c(t)\geqslant (1-\frac{|\xi|}{K_2\gamma})a(t)\label{HB0:ineq1}
\end{equation}
From now on, we  choose: $\lambda=\frac{\gamma K_2}{2\alpha}$. Observe that the constant $K_2$ appearing in the growth condition $\Hd(2)$ can be chosen as small as needed to get: $\lambda \leq\frac{2}{\gamma+2}\alpha<\alpha$. With that choice, we then have:
\begin{equation}
|\xi| = \lambda(\alpha-\lambda)\leqslant \alpha\lambda\leqslant \frac{\gamma K_2}{2},\label{HB0:xi}
\end{equation}
and $2(\lambda-\alpha) \leqslant -\lambda\gamma$. It follows from \eqref{HB0:ineq1} that, for all $t\geq t_1$:
$$a(t)+2\frac{\xi}{\gamma} c(t)\geqslant \frac{1}{2}a(t),$$
and, noticing that  $a(t)+b(t)\geqslant \mathcal E(t)$, we finally get for all $t\geq t_1$:
\begin{eqnarray}
\mathcal E'(t) &\leqslant& -\frac{\lambda\gamma}{2}(a(t) +2b(t)) +\langle g(t),\dot x(t)+\lambda(x(t)-x^*)\rangle\nonumber\\
\mathcal E'(t) &\leqslant &-\frac{\lambda\gamma}{2}\mathcal E(t) +\langle g(t),\dot x(t)+\lambda(x(t)-x^*)\rangle.\label{HB0:diff:ineq1}
\end{eqnarray}
or equivalently:
\begin{equation}
\mathcal G'(t) \leqslant -\frac{\lambda\gamma}{2}\mathcal E(t)\leqslant 0.\label{HB0:diff:ineq2}
\end{equation}

The rest of the proof is quite standard: before integrating the differential inequality \eqref{HB0:diff:ineq1} between $t_1$ and $t$, we first need to control the term $\|\dot x(t) + \lambda(x(t)-x^*)\|$. To that end, observe that, according to \eqref{HB0:diff:ineq2}, the energy $\mathcal G$ is non-increasing, hence: $\forall t\geqslant t_1,~\mathcal G(t) \leq \mathcal G(t_1),$ i.e.:
\begin{eqnarray*}
\forall t\geqslant t_1,~\mathcal E(t) &\leqslant & \mathcal E(t_1) +\int_{t_1}^t \langle g(s),\dot x(s)+\lambda(x(s)-x^*)\rangle ds\\
&\leqslant & \mathcal E(t_1) +\int_{t_1}^t \|g(s)\| \|\dot x(s)+\lambda(x(s)-x^*)\| ds
\end{eqnarray*}
With our choice of parameters, the energy $\mathcal E$ is not a sum of non-negative terms, so that the term $\|\dot x(s)+\lambda(x(s)-x^*)\|$ can not be directly controlled by $\mathcal E(t)$. But according to the growth condition $\Hd(2)$ and to the uniqueness of the minimizer, we have for all $t\geqslant t_1$: $c(t) \leqslant \frac{1}{2K_2}a(t)$. Using \eqref{HB0:xi} and $\gamma\leqslant 2$, we deduce:
$$|\xi| c(t) \leqslant\frac{|\xi|}{2K_2} a(t) \leqslant \frac{\gamma}{4}\leqslant \frac{1}{2},$$
hence: $\forall t\geqslant t_1,~\mathcal E(t) \geqslant \frac{1}{2} a(t) +b(t)\geqslant b(t)$. It follows:
\begin{eqnarray*}
\forall t\geqslant t_1,~b(t) = \|\dot x(t)+\lambda(x(t)-x^*)\|^2 \leqslant \mathcal E(t_1) + \int_{t_1}^t \|g(s)\| \|\dot x(s)+\lambda(x(s)-x^*)\| ds
\end{eqnarray*}
Applying the Gr\"onwall-Bellman Lemma \cite[Lemma A.5]{brezis1973ope}, we obtain:
$$\forall t\geqslant t_1,~\|\dot x(t)+\lambda(x(t)-x^*)\| \leqslant c + \int_{t_1}^t \|g(s)\| ds,$$
where: $c=\sqrt{2\mathcal E(t_1)}$. Since $\int_{t_1}^{+\infty} \|g(s)\|ds <+\infty$ by assumption, we can conclude that:
$$A=\sup_{t\geqslant t_1} \|\dot x(t)+\lambda(x(t) - x^*)\| \leqslant c+\int_{t_1}^{+\infty}\|g(s)\|ds <+\infty.$$
Coming back to \eqref{HB0:diff:ineq1}, we obtain the following differential inequality:
\begin{eqnarray*}
\forall t\geqslant t_1,~\mathcal E'(t)  +\frac{\lambda\gamma}{2}\mathcal E(t)\leqslant   A\|g(t)\|.
\end{eqnarray*}
Integrating between $t_1$ and $t$ , we finally obtain:
\begin{eqnarray*}
\forall t\geqslant t_1,~e^{\frac{\lambda\gamma}{2}t}\mathcal E(t) &\leqslant &e^{\frac{\lambda\gamma}{2}t_1}\mathcal E(t_1) + \int_{t_1}^t e^{\frac{\lambda\gamma}{2}s}\|g(s)\|ds,\\
&\leqslant & e^{\frac{\lambda\gamma}{2}t_1}\mathcal E(t_1) + \int_{t_1}^{+\infty}e^{\frac{\lambda\gamma}{2}s}\|g(s)\|ds <+\infty.
\end{eqnarray*}
Hence: $\mathcal E(t) = \bigO{e^{-\frac{\lambda\gamma}{2}t}}$. Since: $F(x(t))-F^* = a(t) \leqslant 2\mathcal E(t) $ for all $t\geq t_1$, we finally get the expected result.

\paragraph{Case where $\theta\in (0,1)$.} 
Following the strategy proposed by H. Attouch and A. Cabot in the proof of \cite[Theorem 3.12]{attouch2017asymptotic} in the unperturbed case, we choose in time-dependent parameters $\lambda$ and $\xi$:
$$\lambda(t) = \dpy\frac{2\beta(t)}{\gamma+2},~\xi(t) = -\lambda(t)^2,$$
so that the energy $\mathcal E$ can be rewritten as:
\begin{eqnarray*}
\mathcal{E}(t) &=& F(x(t)) - F^* + \frac{1}{2}\|\lambda(x(t) - x^*) + \dot{x}(t)\|^2 + \frac{\xi}{2}\|x(t) - x^*\|^2\\
&=& F(x(t)) - F^* + \frac{1}{2}\| \dot{x}(t)\|^2 + \frac{2\beta(t)}{\gamma+2}\langle x(t) - x^*,\dot x(t)\rangle.
\end{eqnarray*}
The case $\theta=0$ is excluded in the proof detailed hereafter since, as in \cite[Theorem 3.12]{attouch2017asymptotic}, we need that: $\lim_{t\rightarrow +\infty} \beta(t) =0$. 
%
Extending the proof of \cite[Theorem 3.12]{attouch2017asymptotic} to our setting, we obtain the following differential inequality whose proof is detailed in appendix:
\begin{lemma}
Let $\gamma\geq 1$ and $\beta(t) = \frac{\alpha}{t^\theta}$. If $F$ satisfies the hypothesis $\Hu(\gamma)$, then:
\begin{align*}
\mathcal E'(t) +\frac{2\gamma}{\gamma+2}\beta(t)\mathcal E(t)  \leq \frac{2}{\gamma+2}  \left(\dot\beta(t)+\frac{\gamma-2}{\gamma+2}\beta(t)^2\right)&\langle x(t)-x^*,\dot x(t)\rangle\\
& +\left\langle g(t),\dot x(t) + \frac{2}{\gamma+2}\beta(t)(x(t)-x^*)\right\rangle.
\end{align*}\label{lem:ineq:HB}
\end{lemma}
According to Lemma \ref{lem:ineq:HB}  and noticing that for all $t>0$: $\dot\beta(t)+\frac{\gamma-2}{\gamma+2}\beta(t)^2\leqslant 0$, we then obtain:
\begin{eqnarray}
\mathcal G'(t) &\leqslant& -\frac{2\gamma}{\gamma+2}\beta(t)\mathcal E(t)  + \frac{2}{\gamma+2}  \left(\dot\beta(t)+\frac{\gamma-2}{\gamma+2}\beta(t)^2\right)\langle x(t)-x^*,\dot x(t)\rangle,\nonumber\\
&\leqslant &  -\frac{2\gamma}{\gamma+2}\beta(t)\mathcal E(t) + \frac{2}{\gamma+2} \beta(t) \left(\frac{2-\gamma}{\gamma+2}\beta(t)-\frac{\dot\beta(t)}{\beta(t)}\right)\left|\langle x(t)-x^*,\dot x(t)\rangle\right|. \label{HB:ineq1}
\end{eqnarray}
To prove that the energy $\mathcal G$ is non increasing and thus bounded, we have now to control the scalar product $\langle x(t)-x^*,\dot x(t)\rangle$. Assuming that $F$ satisfies the growth condition $\Hd(2)$ and admits a unique minimizer, we first have:
\begin{eqnarray}
\left|\langle x(t)-x^*,\dot x(t)\rangle\right| & \leq & \frac{1}{2}\|x(t)-x^*\|^2 +\frac{1}{2}\|\dot x(t)\|^2
\leqslant  C\left(F(x(t))-F^*+\frac{1}{2}\|\dot x(t)\|^2\right)\label{HB:ineq2}
\end{eqnarray}
where $C=\max(1,\frac{1}{2K_2})$. Since $\dpy\lim_{t\rightarrow +\infty} \beta(t) = 0$, it follows:
\begin{eqnarray*}
\mathcal E(t) &=& F(x(t))-F^* + \frac{1}{2} \|\dot x(t)\|^2 + \frac{2\beta(t)}{\gamma+2} \langle x(t)-x^*,\dot x(t)\rangle\\
&=& F(x(t))-F^* + \frac{1}{2} \|\dot x(t)\|^2 + \smallO{F(x(t))-F^* + \frac{1}{2} \|\dot x(t)\|^2},
\end{eqnarray*}
so that there exists $t_1\geqslant t_0$ such that for all $t\geqslant t_1$, 
\begin{equation}
\mathcal E(t) \geqslant \frac{2}{\gamma+2}\left(F(x(t))-F^* + \frac{1}{2} \|\dot x(t)\|^2\right)(\geqslant 0).\label{HB:ineq3} 
\end{equation}
Combining \eqref{HB:ineq1}, \eqref{HB:ineq2} and \eqref{HB:ineq3}, we then obtain:
$$\mathcal G'(t)  \leqslant -\frac{2\gamma}{\gamma+2}\left(1-C  \left(\frac{2-\gamma}{\gamma+2}\beta(t)-\frac{\dot\beta(t)}{\beta(t)}\right)\right)\beta(t)\mathcal E(t).$$
Observe now that by definition, we have: $\beta(t)^2=\smallO{\beta(t)}$ and $\dot \beta(t)=\smallO{\beta(t)}$. Hence, for any constant $m\in (0,\frac{2\gamma}{\gamma+2})$, there exists $t_2\geq t_1$ such that for all $t\geqslant t_2$:
\begin{equation}
\mathcal G'(t)  \leqslant -m\beta(t)\mathcal E(t)\leqslant 0,\label{diff:ineq1}
\end{equation}
or equivalently:
\begin{equation}
\mathcal E'(t)  + m\beta(t)\mathcal E(t)\leqslant \langle g(t),\dot x(t) + \frac{2}{\gamma+2}\beta(t)(x(t)-x^*)\rangle.\label{diff:ineq2}
\end{equation}

The rest of the proof is quite standard: before integrating the differential inequality \eqref{diff:ineq2} between $t_2$ and $t$, we first need to control the term $\|\dot x(t) + \frac{2}{\gamma+2}\beta(t)(x(t)-x^*)\|$. To that end, observe that, according to \eqref{diff:ineq1}, the energy $\mathcal G$ is non-increasing and that: $\forall t\geqslant t_2,~\mathcal G(t) \leq \mathcal G(t_2),$ i.e.:
\begin{eqnarray*}
\forall t\geqslant t_2,~\mathcal E(t) &\leqslant & \mathcal E(t_2) + \int_{t_2}^t \langle g(s),\dot x(s)+\frac{2\beta(s)}{\gamma+2}(x(s)-x^*)\rangle ds\\
&\leqslant & \mathcal E(t_2) + \int_{t_2}^t \|g(s)\| \|\dot x(s)+\frac{2\beta(s)}{\gamma+2}(x(s)-x^*)\|ds.
\end{eqnarray*}
Moreover, with our choice of parameters, the energy $\mathcal E(t)$ is not a sum of non negative terms:
\begin{equation}
\mathcal{E}(t) =F(x(t))-F^* + \frac{1}{2}\|\dot x(t)+\frac{2\beta(t)}{\gamma+2}(x(t) - x^*)\|^2 - \frac{2\beta(t)^2}{(\gamma+2)^2}\|x(t) - x^*\|^2,
\end{equation}
so that the term $\|\dot x(t)+\frac{2\beta(t)}{\gamma+2}(x(t) - x^*)\|^2$ can not be directly controlled by the energy $\mathcal E(t)$. But, according to the growth condition $\Hd(2)$ combined with the uniqueness of the minimizer, we have:
\begin{eqnarray*}
F(x(t)) - F^* - \frac{2\beta(t)^2}{(\gamma+2)^2}\|x(t) - x^*\|^2 &\geqslant& \left(1-\frac{2\beta(t)^2}{K_2(\gamma+2)^2}\right)(F(x(t))-F^*)\\
&\geqslant& \frac{1}{2}(F(x(t))-F^*).
\end{eqnarray*}
for $t$ large enough, and: $\mathcal E(t) \geqslant  \frac{1}{2}(F(x(t))-F^*) + \|\dot x(t)+\frac{2\beta(t)}{\gamma+2}(x(t) - x^*)\|^2$. Hence:
\begin{eqnarray*}
\frac{1}{2}\|\dot x(t)+\frac{2\beta(t)}{\gamma+2}(x(t) - x^*)\|^2 &\leqslant& \mathcal E(t)\leqslant \mathcal E(t_2) + \int_{t_2}^t \langle g(s),\dot x(s)+\frac{2\beta(s)}{\gamma+2}(x(s)-x^*)\rangle ds\\
&\leqslant &  \frac{1}{2}c^2 + \int_{t_2}^t \|g(s)\|\|\dot x(s)+\frac{2\beta(s)}{\gamma+2}(x(s)-x^*)\| ds
\end{eqnarray*}
where: $c=\sqrt{2\mathcal E(t_2)}$. Applying the Gr\"onwall-Bellman Lemma \cite[Lemma A.5]{brezis1973ope}, we obtain:
$$\|\dot x(t)+\frac{2\beta(t)}{\gamma+2}(x(t) - x^*)\| \leq c +\int_{t_2}^t\|g(s)\|ds.$$
Assuming that $\int_{t_2}^{+\infty}\|g(s)\|ds<+\infty$, we can so conclude that:
$$A=\sup_{t\geqslant t_2} \|\dot x(t)+\frac{2\beta(t)}{\gamma+2}(x(t) - x^*)\| \leqslant c+\int_{t_2}^{+\infty}\|g(s)\|ds <+\infty.$$
Coming back to \eqref{diff:ineq2}, we obtain the following differential inequality:
\begin{eqnarray*}
\forall t\geqslant t_2,~\mathcal E'(t)  + m\beta(t)\mathcal E(t)&\leqslant & \|g(t)\|\|\dot x(t) + \frac{2\beta(t)}{\gamma+2}(x(t)-x^*)\|\\
&\leqslant & A\|g(t)\|.
\end{eqnarray*}
Integrating between $t_2$ and $t$ and stating: $\Gamma (t) = \dpy\int_{t_2}^{t} \beta(s)ds$, we finally obtain: for all $t\geqslant t_2$,
\begin{eqnarray*}
e^{m\Gamma(t)}\mathcal E(t) &\leqslant &e^{m\Gamma(t_2)}\mathcal E(t_2) + A\int_{t_2}^t e^{m\Gamma(s)}\|g(s)\|ds,\\
&\leqslant & e^{m\Gamma(t_2)}\mathcal E(t_2) + A\int_{t_2}^{+\infty}e^{m\Gamma(s)}\|g(s)\|ds =B<+\infty.
\end{eqnarray*}
Hence: $\mathcal E(t) = \bigO{e^{-m\Gamma(t)}}$. Since: $F(x(t))-F^* = a(t) \leqslant 2\mathcal E(t) $ for all $t\geq t_2$, we finally get the expected result and using $\Hd(2)$, the other estimates follow directly.

\subsection{Proof of Theorem \ref{th:unperturbed:flat case}}
The proof of Theorem \ref{th:unperturbed:flat case} relies on almost the same energy as that used in \cite{aujol2018optimal}:
\begin{equation}
\mathcal{E}(t) = t^{2}(F(x(t)) - F^*) + \frac{1}{2}\|\lambda(x(t) - x^*) + t\dot{x}(t)\|^2 + \frac{\xi(t)}{2} \|x(t) - x^*\|^2
\end{equation}
where $\lambda$ is a non-negative real constant as in \cite{aujol2018optimal} and $\xi(.)$ is here a real-valued function. 
Noting:
\begin{equation}
    \begin{aligned}
    a(t) & = t(F(x(t)) - F^*)\\
    b(t) & = \frac{1}{2t}\|\lambda(x(t) - x^*) + t\dot{x}(t)\|^2\\
    c(t) & = \frac{1}{2t}\|x(t) - x^*\|^2
    \end{aligned}
\end{equation}
we have: $\mathcal{E}(t) = t(a(t) + b(t) + \xi(t)c(t)) $. We also define:
\begin{equation}
    \mathcal{H}(t) = t^p\mathcal{E}(t)
\end{equation}
The proofs of our theorems rely on the following lemma whose proof is detailed in appendix \ref{Appendix HB}:

\begin{lemma} \label{lemma HB}
Let $\gamma_1 \geqslant 1$. If $F$ satisfies the hypothesis $\mathbf{H}_1(\gamma_1)$ and if $\xi(t) = \lambda(\lambda + 1 - \alpha t^{1-\theta})$, then:
\begin{equation} \label{thm 4 central in proof}
\begin{aligned}
\mathcal{H}'(t) & \leqslant t^p\Bigg((2 + p - \gamma_1 \lambda)a(t) + (2\lambda +2 +p -2\alpha t^{1-\theta})b(t)+ \lambda((\lambda+1)(p-2\lambda)-\alpha(p+1-\theta-2\lambda)t^{1-\theta}) c(t)\Bigg)
\end{aligned}
\end{equation}
\end{lemma}

Note $r=\frac{1+\theta}{2}$. Taking $\lambda = \frac{2r}{\gamma_1 - 2}$ and $p =p_1+2(r-1)$ with $p_1=\frac{4r}{\gamma_1 - 2}$, we obtain:
\begin{equation}
\begin{aligned}
\mathcal{H}'(t) \leqslant t^{p}\Bigg(2\Big(\frac{\gamma_1+2}{\gamma_1-2}r - \alpha t^{-2(r-1)}\Big)b(t) + 2\lambda(\lambda+1)(r-1) c(t)\Bigg)
\end{aligned}
\end{equation}
Since $r\leqslant 1$, $2\lambda(\lambda+1)(r-1) \leqslant 0$, hence:
\begin{equation} \label{equation c for r=1}
\begin{aligned}
\mathcal{H}'(t) \leqslant 2t^{p_1} \Big(\frac{\gamma_1+2}{\gamma_1-2}rt^{2(r-1)} - \alpha\Big)b(t)
\end{aligned}
\end{equation}
\begin{enumerate}
\item If $\theta<1$, it exists $t_1$ depending only on $\gamma_1$, $\alpha$ and $\theta$ such that $\mathcal{H}'(t)\leqslant 0$ for all $t\geqslant t_1$. 
\item If $\theta=1$ and $\alpha\geqslant \frac{\gamma_1+2}{\gamma_1-2}$, $\mathcal{H}'$ is non positive for all $t\geqslant t_0$.
\end{enumerate}

We will now use similar reasoning as in \cite{aujol2018optimal} to prove the results of our theorem. Since $\mathcal{H}'(t) \leqslant 0$, for any choice of $x^*$
in the set of minimizers $X^*$, the function $\mathcal{H}$ is bounded above and since the set of minimizers is bounded because $F$ is coercive, there exists $A$ > 0 and $t_0$ such that for all choices of $x^*$
in $X^*$:
\begin{equation}
\mathcal{H}(t_0) \leqslant A
\end{equation}
Hence for all $x^* \in X^*$ and $t \geqslant t_0$, $\mathcal{H}(t) \leqslant A$. Hence
\begin{equation}
t^{\frac{2r\gamma_1}{\gamma_1 - 2}}(F(x(t)) - F^*) \leqslant \frac{|\xi(t)|}{2}t^{\frac{4r}{\gamma_1 - 2}+2(r-1)}\|x(t) - x^*\|^2 + A.
\end{equation}
We have $\xi(t) = \frac{2r}{\gamma_1 - 2}( \frac{2r}{\gamma_1 - 2}+1 - \alpha t^{-2(r-1)})$. Hence:
\begin{itemize}
    \item If $\theta =1$, then $r=1$. Since $\alpha \geqslant \frac{\gamma_1 + 2}{\gamma_1 - 2} \geqslant \frac{\gamma_1}{\gamma_1 - 2}$, we have:  
    $$|\xi(t)| =  \frac{2}{\gamma_1-2}(\alpha -\frac{\gamma_1}{\gamma_1-2})\leqslant \frac{2\alpha}{\gamma_1-2}.$$
    \item If $\theta\in (0,1)$, then:
    $$|\xi(t)|t^{2(r-1)} =\frac{2r}{\gamma_1-2}(\alpha -(\frac{2r}{\gamma_1-2}+1)t^{2(r-1)}\leqslant \frac{2r\alpha}{\gamma_1-2}.$$
\end{itemize}
Hence, in both cases:
\begin{equation}\label{Bornexi}
|\xi(t)|t^{2(r-1)} \leqslant \frac{2r\alpha}{\gamma_1 - 2}.
\end{equation}
Therefore:
\begin{equation}
t^{\frac{2r\gamma_1}{\gamma_1 - 2}}(F(x(t)) - F^*) \leqslant \frac{r\alpha}{\gamma_1 - 2}t^{\frac{4r}{\gamma_1 - 2}}\|x(t) - x^*\|^2 + A
\end{equation}
And since this is verified for all $x^* \in X^*$:
\begin{equation}
t^{\frac{2r\gamma_1}{\gamma_1 - 2}}(F(x(t)) - F^*) \leqslant \frac{r\alpha}{\gamma_1 - 2}t^{\frac{4r}{\gamma_1 - 2}}d(x(t), X^*)^2 + A
\end{equation}
We set $v(t) = t^{\frac{4r}{\gamma_2 - 2}} d(x(t), X^*)^2$
Then
\begin{equation} \label{HBF proof closing}
t^{\frac{2r\gamma_1}{\gamma_1 - 2}}(F(x(t)) - F^*) \leqslant \frac{r\alpha}{\gamma_1 - 2}t^{\frac{4r}{\gamma_1 - 2} - \frac{4r}{\gamma_2 - 2}}v(t) + A
\end{equation}
Since $F$ satisfies $\mathbf{H}_2(\gamma_2)$, there exists K > 0 such that
\begin{equation}
K(t^{-\frac{4r}{\gamma_2 - 2}}v(t))^{\frac{\gamma_2}{2}} \leqslant F(x(t)) - F^*
\end{equation}
i.e
\begin{equation}
Kv(t)^{\frac{\gamma_2}{2}}t^{\frac{-2r\gamma_2}{\gamma_2 - 2}} \leqslant F(x(t)) - F^*
\end{equation}
Hence
\begin{equation}
Kt^{\frac{2r\gamma_1}{\gamma_1 - 2}} t^{- \frac{2r\gamma_2}{\gamma_2 - 2}}v(t)^{\frac{\gamma_2}{2}} \leqslant t^{\frac{2r\gamma_1}{\gamma_1 - 2}}(F(x(t)) - F^*)
\end{equation}
Back to (\ref{HBF proof closing}), this yields:
\begin{equation}
Kt^{\frac{2r\gamma_1}{\gamma_1 - 2}} t^{- \frac{2r\gamma_2}{\gamma_2 - 2}}v(t)^{\frac{\gamma_2}{2}} \leqslant \frac{r\alpha}{\gamma_1 - 2} t^{\frac{4r}{\gamma_1 - 2} - \frac{4r}{\gamma_2 - 2}}v(t) + A
\end{equation}
Hence
\begin{equation}
Kv(t)^{\frac{\gamma_2}{2}} \leqslant \frac{rc}{\gamma - 2}v(t) + At^{\frac{4r}{\gamma_2 - 2} - \frac{4r}{\gamma_1 - 2}}
\end{equation}
Which, since $\gamma_1 \leqslant \gamma_2$, means that $v$ is bounded. Therefore, from (\ref{HBF proof closing}) we deduce that there exists $B > 0$ such that:
\begin{equation}
F(x(t)) - F^* \leqslant Bt^{\frac{-2r\gamma_2}{\gamma_2 - 2}} + At^{\frac{-2r\gamma_1}{\gamma_1 - 2}}
\end{equation}
Since $\gamma_1 \leqslant \gamma_2$, we have $\frac{-r2\gamma_2}{\gamma_2 - 2} \geqslant \frac{-2r\gamma_1}{\gamma_1 - 2}$. Hence $F(x(t)) - F^* = O(t^{-\frac{2r\gamma_2}{\gamma_2 - 2}})$.


\subsection{Proof of Theorem \ref{th:perturbed:flat case}}
The proof is inspired by the one of \textbf{Theorem} \ref{th:unperturbed:flat case} where an additional term including the noise is considered. 

First, we set parameters $r=\frac{1+\theta}{2}$, $\lambda=\frac{2r}{\gamma_1-2}$ and $p=p_1+2(r-1)$, with $p_1=2\lambda$ and functions $\xi$, $\mathcal{E}$ and $\mathcal{H}$ exactly as in \textbf{Theorem} \ref{th:unperturbed:flat case}. 
In addition to the energy functions $\mathcal{E}$ and $\mathcal{H}$ we define
\begin{equation}
\mathcal{G}(t) = \mathcal{H}(t) + \int_{t}^T  \langle \lambda(x(s) - x^*) + s\dot{x}(s) , s^{p+1}g(s)\rangle ds
\end{equation}
We refer the reader to Appendix \eqref{Appendix HB} for a detailed calculation of the following bound on the derivative $\mathcal{G}'$  
\begin{equation}
\mathcal{G}'(t) \leqslant t^p\Bigg((2 + p - \gamma_1 \lambda)a(t) + (2\lambda +2 +p -2\alpha t^{1-\theta})b(t)+ \lambda((\lambda+1)(p-2\lambda)-\alpha(p+1-\theta-2\lambda)t^{1-\theta}) c(t)\Bigg)
\end{equation}
Actually, the integral term 
$\mathcal{G}(t) -\mathcal{H}(t)$ is computed such that the bound on $\mathcal{G}'(t)$ is equal to the bound \eqref{thm 4 central in proof} proposed in  \textbf{Theorem} \ref{th:unperturbed:flat case} for $\mathcal{H}'$. Thus we deduce once again that under the hypotheses of the Theorem \ref{th:perturbed:flat case}, it exists $t_1\geqslant 1$ such that the function $\mathcal{G}$ is non increasing for $t\geqslant t_1$.

We will need now the following direct lemma
\begin{lemma}\label{LemmeBornec}
If $F$ satisfies the growth condition $\mathbf{H}_2(\gamma_2)$ with $\gamma_2>2$ :
\begin{equation}
K\|x-x^*\|^{\gamma_2}\leqslant F(x)-F^*
\end{equation}
then defining $p_2=\frac{4r}{\gamma_2-2}$, with notations defined in \eqref{eqdefabc} we have 
\begin{equation}
    t^{p_2+1}c(t)\leqslant \frac{K^{-\frac{2}{\gamma_2}}}{2}\left(t^{p_2+2r-1}a(t)\right)^{\frac{2}{\gamma_2}}
\end{equation}
It follows that for any $m\in\mathbb{R}$, it exists $M\in \mathbb{R}$ such that for any $t\geqslant t_0$
\begin{equation}
    mt^{p_2+1}c(t)-t^{p_2+2r-1}a(t)\leqslant M.
\end{equation}
\end{lemma}

Since for all $t \geqslant t_1$, $\mathcal{G}(t) \leqslant \mathcal{G}(t_1)$. Then:
\begin{equation}
\mathcal{H}(t) \leqslant \mathcal{H}(t_1) + \int_{t_1}^t \langle (\lambda (x(s) - x^*) + s \dot{x}(s)), s^{p+1} g(s) \rangle ds.
\end{equation}
Hence, we have:
\begin{equation}
t^{p+1}a(t) + t^{p+1}b(t) \leqslant \mathcal{H}(t_1) + |\xi(t)|t^{p+1}c(t) + \int_{t_1}^t \langle  \lambda (x(s) - x^*) + s \dot{x}(s), s^{p+1} g(s) \rangle ds.
\end{equation}
Using the fact that the $|\xi(t)|t^{2(r-1)}$ is uniformly bounded, see \eqref{Bornexi} we get 
\begin{equation}
t^{p+1}a(t) + t^{p+1}b(t) \leqslant \mathcal{H}(t_1) + \frac{2r\alpha}{\gamma_1 - 2}t^{p_1+1}c(t) + \int_{t_1}^t \| \lambda (x(s) - x^*) + s\dot{x}(s)\|_2 \|s^{p+1} g(s)\|_2 ds 
\end{equation}
Let's define $p_2=\frac{4r}{\gamma_2-2}$. Since $\gamma_1\leqslant \gamma_2$ we have $p_2\leqslant p_1$ and for any $t\geqslant t_1$, 
$t^{p_1-p_2}\geqslant 1$. Dividing the previous inequality by 
$t^{p_1-p_2}$ we get for any $t\geqslant t_1$ 
\begin{equation}
t^{p_2+2r-1}(a(t) + b(t)) \leqslant \mathcal{H}(t_1) + \frac{2r\alpha}{\gamma_1 - 2}t^{p_2+1}c(t) + t^{p_2-p_1}\int_{t_1}^t \| \lambda (x(s) - x^*) + s \dot{x}(s)\|_2 \|s^{p+1} g(s)\|_2 ds 
\end{equation}
which implies 
\begin{equation}\label{eqp2}
t^{p_2+2r-1}(a(t) + b(t)) \leqslant \mathcal{H}(t_1) + \frac{2r\alpha}{\gamma_1 - 2}t^{p_2+1}c(t) + \int_{t_1}^t \| \lambda (x(s) - x^*) + s \dot{x}(s)\|_2 \|s^{p_2+2r-1} g(s)\|_2 ds 
\end{equation}
Since $F$ satisfies the growth condition $\mathbf{H}(\gamma_2)$ we can apply Lemma \ref{LemmeBornec} and deduce it exists $M\in\mathbb{R}$ such that for $t\geqslant t_1$
\begin{equation}
t^{p_2+2r-1}b(t) \leqslant M+ \int_{t_1}^t \| (\lambda(x(s) - x^*) + s \dot{x}(s)\|_2 \|s^{p_2+2r-1} g(s)\| \rangle ds 
\end{equation}
which implies using the definition of $b(t)$ given in \eqref{eqdefabc} 
\begin{equation}
\frac{1}{2}t^{p_2+2(r-1)}\|\lambda(x(t)-x^*)+t\dot x(t)\|_2^2 \leqslant M+ \int_{t_1}^t \|s^{\frac{p_2}{2}+(r-1)} (\lambda (x(s) - x^*) + s \dot{x}(s))\|_2 \|s^{\frac{p_2}{2}+r} g(s)\|_2ds 
\end{equation}
Applying the Gr\"onwall Bellman Lemma it follows that 
\begin{equation}
t^{\frac{p_2}{2}+r-1}\|\lambda(x(t)-x^*)+t\dot x(t)\|_2 \leqslant \sqrt{2M}+ \int_{t_1}^t s^{\frac{p_2}{2}+r} \|g(s)\|_2 ds 
\end{equation}
Under the hypotheses of the Theorem, the right member of the inequality is uniformly bounded relatively to 
$t$. It follows that it exists $M_1\geqslant 0$ such that for any $t\geqslant t_1$
\begin{equation}
\|\lambda(x(t)-x^*)+t\dot x(t)\|_2 \leqslant M_1t^{-\frac{p_2}{2}-r+1}\label{boundflat}
\end{equation}
and thus that it exists $M_2>0$ such that for any $t\geqslant t_1$
\begin{equation}
\int_{t_1}^t |\langle  (\lambda (x(s) - x^*) + s \dot{x}(s)), s^{p_2+2r-1} g(s) \rangle| ds\leqslant M_1 \int_{t_1}^t s^{\frac{p_2}{2}+r}\|g(s)\|ds\leq M_2
\end{equation}

Combining this inequality with \eqref{eqp2} if follows that for any $t\geqslant t_1$
\begin{equation}
t^{p_2+2r-1}a(t) \leqslant \mathcal{H}(t_1) + \frac{2r\alpha}{\gamma_1 - 2}t^{p_2+1}c(t) +M_2
\end{equation}
Using once again Lemma \ref{LemmeBornec} we deduce it exists $M_3$ and $M_4$ such that for any $t\geqslant t_1$
\begin{equation}
t^{p_2+2r-1}a(t) \leqslant M_3+ M_4(t^{(p_2+2r-1)}a(t))^{\frac{\gamma_2}{2}}
\end{equation}
which implies that it exists $M_5$ such that any $t\geqslant 1$
\begin{equation}
a(t) \leqslant M_5t^{-(p_2+2r-1)}
\end{equation}
that is 
\begin{equation}
F(x(t))-F^*\leqslant M_5t^{-(p_2+2r)}=M_5t^{-\frac{2r\gamma_2}{\gamma_2-2}}
\end{equation}
Which is the desired result.

\subsection{Proof of Corollary \ref{cor:flat}}
In this paragraph we detail the proof of Corollary \ref{cor:flat}  in which it is stated that, in the flat case, the trajectory of any solution $x$ of the ODE (\ref{ODE:HB}) is finite.

From the proof of Theorem \ref{th:perturbed:flat case} (see \eqref{boundflat}), there exists $A_1 > 0$ and $t_1\geqslant t_0$ such that, for all $t\geqslant t_1$:
\begin{equation}
\|\lambda(x(t) - x^*) + t\dot{x}(t)\| \leqslant A_1 t^{-\frac{r\gamma}{\gamma - 2}+1}
\end{equation}
Combining the growth condition $\mathbf{H}_2(\gamma)$ and the conclusion of Theorem \ref{th:perturbed:flat case}, we have that there exists $t_2\geqslant t_1$ such that:
\begin{eqnarray*}
\forall t\geqslant t_2,~\|x(t) - x^*\| &\leqslant & K^{-\frac{1}{\gamma}}(F(x(t)) - F^*)^{\frac{1}{\gamma}}\\
&\leqslant & A_2 t^{-\frac{2r}{\gamma - 2}}
\end{eqnarray*}
It follows that for all $t\geqslant t_2$:
\begin{eqnarray*} \label{central cor 1}
t\|\dot{x}(t)\| &\leqslant& \lambda  \|x(t) - x^*\| + \|\lambda(x(t) - x^*) + t\dot{x}(t)\|\\
&\leqslant&  A_1 t^{-\frac{r\gamma}{\gamma - 2}+1} + \lambda A_2 t^{-\frac{2r}{\gamma - 2}}\\
&\leqslant &  t^{-\frac{r\gamma}{\gamma - 2}+1}(A_1 + \lambda A_2 t^{r-1})
\end{eqnarray*}
Hence:
$$\|\dot{x}(t)\| \leqslant t^{-\frac{r\gamma}{\gamma - 2}}(A_1 + \lambda A_2 t^{r-1}).$$
Noticing that $r-1\leqslant 0$, we finally get:
\begin{equation*}
\|\dot{x}(t)\| =\bigO{t^{-\frac{r\gamma}{\gamma - 2}}}.
\end{equation*}
Which means that $\|\dot{x}(t)\|$ is integrable and that the trajectory is finite.

\section*{Acknowledgement}
This study has been carried out with financial support from the French state, managed by the French National Research Agency (ANR GOTMI - ANR-16-VCE33-0010-01) and partially supported by ANR-11-LABX-0040-CIMI within the program
ANR-11-IDEX-0002-02. 

\appendix

\section{Appendix}\label{sec_appendix}
\subsection{Proof of Lemma \ref{lemma Nesterov}: differentiating the energy function for Theorem \ref{th:Nesterov:sharp2}} \label{Appendix Nesterov}
\begin{enumerate}
\item \textbf{Differentiating} $\mathbf{\mathcal{E}}$.

Consider the energy $\mathcal E$ defined as follows by:
\begin{eqnarray*}
\mathcal E(t) &=& t^2(F(x(t))-F^*) +\frac{1}{2}\|\lambda(x(t)-x^*)+t\dot x(t)\|^2 +\frac{\xi}{2}\|x(t)-x^*\|^2\\
&=& t(a(t) + b(t) +c(t)) 
\end{eqnarray*}
where:
\begin{equation}
a(t) = t(F(x(t)) - F^*),~b(t) = \frac{1}{2t}\|\lambda(x(t) - x^*) + t\dot{x}(t)\|^2,~c(t) = \frac{1}{2t}\|x(t) - x^*\|^2.
\end{equation}
We then have:
\begin{equation}
\begin{aligned}
\mathcal{E}'(t)  = 2t(F(x(t)) - F^*) &+ t^2 \langle \nabla F(x(t)), \dot{x}(t)\rangle + \xi \langle \dot{x}(t), x(t) - x^* \rangle  \\
&+ \langle \lambda(x(t) - x^*) + t \dot{x}(t), (\lambda+1) \dot{x}(t) + t \ddot{x}(t) \rangle.
\end{aligned}
\end{equation}
Since $x$ is a solution of the ODE (\ref{ODE:Nesterov}), we have:
\begin{equation}
\begin{aligned}
(\lambda + 1)\dot{x}(t) + t\ddot{x}(t) &  = (\lambda + 1) \dot{x}(t) - \alpha \dot{x}(t) - t\nabla F(x(t)) + tg(t) \\
& = (\lambda + 1 - \alpha)\dot{x}(t) - t \nabla F(x(t)) + t g(t). 
\end{aligned}
\end{equation}
Hence:
\begin{equation}
\begin{aligned}
\mathcal{E}'(t)  = 2a(t) - \lambda t \langle \nabla F(x(t)), x(t) - x^* \rangle &+ (\xi + \lambda(\lambda + 1 - \alpha)) \langle \dot{x}(t), x(t) - x^*\rangle +  t(\lambda + 1 - \alpha) \|\dot{x}(t)\|^2 \\
& + \langle tg(t), \lambda (x(t) - x^*) + t\dot{x}(t) \rangle.
\end{aligned}
\end{equation}
 
Noticing that:
\begin{equation}
\begin{aligned}
\frac{1}{t} \|\lambda(x(t) - x^*) + t \dot{x}(t) \|^2 = t\|\dot{x}(t)\|^2 + 2\lambda \langle \dot{x}(t), x(t) - x^* \rangle + \frac{\lambda^2}{t} \|x(t) - x^*\|^2,
\end{aligned}
\end{equation}
we deduce that
\begin{eqnarray*}
\mathcal{E}'(t) &=& 2a(t) - \lambda t \langle \nabla F(x(t)), x(t) - x^* \rangle + (\xi - \lambda(\lambda + 1 - \alpha)) \langle \dot{x}(t), x(t) - x^*\rangle \\
&&  + \frac{\lambda + 1 - \alpha}{t} \|\lambda(x(t) - x^*) + t \dot{x}(t) \|^2 - \frac{\lambda^2(\lambda + 1 - \alpha)}{t} \|x(t) - x^*\|^2 \\
&&  + \langle tg(t), \lambda (x(t) - x^*) + t\dot{x}(t) \rangle\\
&=& 2a(t) - \lambda t \langle \nabla F(x(t)), x(t) - x^* \rangle + (\xi - \lambda(\lambda + 1 - \alpha)) \langle \dot{x}(t), x(t) - x^*\rangle \\
& & + 2(\lambda + 1 - \alpha)b(t) - 2\lambda^2(\lambda + 1 - \alpha)c(t) + \langle tg(t), \lambda (x(t) - x^*) + t\dot{x}(t) \rangle.
\end{eqnarray*}
Choosing now $\underline{\xi = \lambda(\lambda + 1 - \alpha)}$, we get:
\begin{equation*}
\begin{aligned}
\mathcal{E}'(t)  = 2a(t) &- \lambda t \langle \nabla F(x(t)), x(t) - x^* \rangle + 2(\lambda + 1 - \alpha)b(t) - 2\lambda^2(\lambda + 1 - \alpha)c(t) \\
& + \langle tg(t), \lambda (x(t) - x^*) + t\dot{x}(t) \rangle
\end{aligned}
\end{equation*}
Since $F$ satisfies $\mathbf{H_1(\gamma)}$:
\begin{equation} \label{dE/dt}
\mathcal{E}'(t) \leqslant (2 - \lambda \gamma)a(t) + 2(\lambda + 1 - \alpha)b(t) - 2\lambda^2(\lambda + 1 - \alpha)c(t) + \langle tg(t), \lambda (x(t) - x^*) + t\dot{x}(t) \rangle
\end{equation}
\item \textbf{Differentiating} $ \mathcal{H} $ and $\mathcal{G}$

Recall now that:
\begin{equation*}
\mathcal{G}(t) = \mathcal{H}(t) + \int_t^T s^p\langle  (\lambda (x(s) - x^*) + s \dot{x}(s)), s g(s) \rangle ds
\end{equation*}
where: $\mathcal{H}(t) =t^p \mathcal{E}(t) $. Since $\mathcal{E}(t) = t(a(t) + b(t) + \xi c(t))$, we deduce from (\ref{dE/dt}) that:
\begin{eqnarray*}
\mathcal{H}'(t) &=& t^{p-1}(p\mathcal{E}(t)+t\mathcal E'(t))\\
&\leqslant & t^p\left((2 - \gamma \lambda + p) a(t) + (2 \lambda +2 - 2 \alpha + p)b(t) + \lambda (\lambda + 1 - \alpha)(-2\lambda + p)c(t)\right) \\
&& + \langle t^{p+1}g(t), \lambda (x(t) - x^*) + t\dot{x}(t) \rangle
\end{eqnarray*}
Hence the expected inequality:
\begin{eqnarray*}
\mathcal{G}'(t) &=& \mathcal H'(t) - t^p\langle  (\lambda (x(t) - x^*) + t \dot{x}(t)), t g(t) \rangle\\
&\leqslant &t^p((2 - \gamma \lambda + p) a(t) + (2 \lambda +2 - 2 \alpha + p)b(t) + \lambda (\lambda + 1 - \alpha)( p-2\lambda)c(t))
\end{eqnarray*}

\end{enumerate}

\subsection{Proof of Lemma \ref{lem:HB:sharp}: differentiating the energy function of Theorem \ref{th:HB:sharp}}
Consider the energy $\mathcal E$ defined as follows by:
\begin{eqnarray*}
\mathcal E(t) &=& F(x(t))-F^* +\frac{1}{2}\|\lambda(x(t)-x^*)+\dot x(t)\|^2 +\frac{\xi}{2}\|x(t)-x^*\|^2\\
&=& a(t) + b(t) +c(t)
\end{eqnarray*}
where:
\begin{equation}
a(t) = F(x(t)) - F^*,~b(t) = \frac{1}{2}\|\lambda(x(t) - x^*) + t\dot{x}(t)\|^2,~c(t) = \frac{1}{2}\|x(t) - x^*\|^2.
\end{equation}

We then have:
\begin{equation}
\mathcal{E}'(t)  =  \langle \nabla F(x(t)), \dot{x}(t)\rangle + \langle \lambda(x(t) - x^*) +  \dot{x}(t), \lambda\dot{x}(t) +  \ddot{x}(t) \rangle+ \xi \langle x(t) - x^*,\dot{x}(t)\rangle.
\end{equation}
Since $x$ is a solution of the ODE (\ref{ODE:HB}), we have:
\begin{equation*}
\lambda\dot{x}(t) + \ddot{x}(t)   = (\lambda -\alpha) \dot{x}(t) - \nabla F(x(t)) + g(t) 
\end{equation*}
Hence:
\begin{eqnarray*}
\mathcal{E}'(t)  &=& - \lambda \langle \nabla F(x(t)), x(t) - x^* \rangle + (\xi + \lambda(\lambda - \alpha))\langle \dot{x}(t), x(t) - x^*\rangle  +(\lambda  - \alpha) \|\dot{x}(t)\|^2 \\
&& + \langle g(t), \lambda (x(t) - x^*) + \dot{x}(t) \rangle.
\end{eqnarray*}
Noticing that:
\begin{equation}
\begin{aligned}
 \|\lambda(x(t) - x^*) + \dot{x}(t) \|^2 = \|\dot{x}(t)\|^2 + 2\lambda \langle \dot{x}(t), x(t) - x^* \rangle + \lambda^2 \|x(t) - x^*\|^2,
\end{aligned}
\end{equation}
we deduce that
\begin{eqnarray*}
\mathcal{E}'(t) &=&  - \lambda  \langle \nabla F(x(t)), x(t) - x^* \rangle + (\lambda - \alpha)\|\lambda(x(t) - x^*) + \dot{x}(t) \|^2+ (\xi - \lambda(\lambda  - \alpha)) \langle \dot{x}(t), x(t) - x^*\rangle \\
&&   - \lambda^2(\lambda - \alpha) \|x(t) - x^*\|^2  + \langle g(t), \lambda (x(t) - x^*) + \dot{x}(t) \rangle\\
&=& - \lambda  \langle \nabla F(x(t)), x(t) - x^* \rangle + 2(\lambda - \alpha)b(t)+ (\xi - \lambda(\lambda  - \alpha)) \langle \dot{x}(t), x(t) - x^*\rangle \\
&&- 2\lambda^2(\lambda  - \alpha)c(t) + \langle g(t), \lambda (x(t) - x^*) + \dot{x}(t) \rangle.
\end{eqnarray*}
Choosing now $\underline{\xi = \lambda(\lambda  - \alpha)}$, we get:
\begin{eqnarray*}
\mathcal{E}'(t) &=&- \lambda  \langle \nabla F(x(t)), x(t) - x^* \rangle + 2(\lambda - \alpha)b(t)- 2\lambda^2(\lambda  - \alpha)c(t) + \langle g(t), \lambda (x(t) - x^*) + \dot{x}(t) \rangle.
\end{eqnarray*}
Since $F$ satisfies $\Hu(\gamma)$, we finally get:
\begin{eqnarray*}
\mathcal{E}'(t) &\leqslant &- \lambda \gamma a(t) + 2(\lambda - \alpha)b(t)- 2\lambda^2(\lambda  - \alpha)c(t) + \langle g(t), \lambda (x(t) - x^*) + \dot{x}(t) \rangle.
\end{eqnarray*}

\subsection{Proof of Lemma \ref{lem:ineq:HB}: differentiating the energy function of Theorem \ref{th:HB:sharp} for the classical heavy ball system}
Consider the energy:
\begin{eqnarray*}
\mathcal{E}(t) &=& F(x(t)) - F^* + \frac{1}{2}\|\lambda(t)(x(t) - x^*) + \dot{x}(t)\|^2 + \frac{\xi(t)}{2}\|x(t) - x^*\|^2\\
&=& F(x(t)) - F^* + \frac{1}{2}\| \dot{x}(t)\|^2 + \frac{2\beta(t)}{\gamma+2}\langle x(t) - x^*,\dot x(t)\rangle.
\end{eqnarray*}
where: $\lambda(t) = \frac{2\beta}{\gamma+2}$ and $\xi(t) = -\lambda(t)^2$. Then:
\begin{eqnarray*}
\mathcal E'(t) &=& \langle \nabla F(x(t)),\dot x(t)\rangle +\langle \ddot x(t) +\frac{2\beta(t)}{\gamma+2} \dot x(t),\dot x(t)\rangle + \frac{2\beta(t)}{\gamma+2} \langle x(t)-x^*,\ddot x(t)\rangle + \frac{2\dot\beta(t)}{\gamma+2} \langle x(t)-x^*,\dot x(t)\rangle.
\end{eqnarray*}
Since $x$ satisfies the ODE \eqref{ODE:HB} when $\theta\in (0,1)$, we have:
$$\ddot x(t) = g(t)-\beta(t)\dot x(t) -\nabla F(x(t)),$$
so that:
\begin{eqnarray*}
\mathcal{E}'(t) &= & -\frac{2\beta(t)}{\gamma+2} \langle \nabla F(x(t)),x(t)-x^*\rangle -\frac{\gamma}{\gamma+2}\beta(t)\|\dot x(t)\|^2 + \frac{2}{\gamma+2}\left(\dot\beta(t) -\beta(t)^2\right)\langle x(t)-x^*,\dot x(t)\rangle\\
&&+ \langle g(t),\dot x(t) +\frac{2\beta(t)}{2+\gamma}(x(t)-x^*)\rangle\\
&\leqslant & -\frac{2\gamma}{\gamma+2}\beta(t)\left[ F(x(t))-F^* +\frac{1}{2}\|\dot x(t)\|^2\right] + \frac{2\beta(t)}{\gamma+2}\left(\frac{\dot\beta(t)}{\beta(t)} -\beta(t)\right)\langle x(t)-x^*,\dot x(t)\rangle\\
&&+ \langle g(t),\dot x(t) +\frac{2\beta(t)}{2+\gamma}(x(t)-x^*)\rangle
\end{eqnarray*}
assuming that $F$ satisfies $\Hu(\gamma)$. Noticing that by definition of the energy $\mathcal E(t)$:
$$F(x(t))-F^* +\frac{1}{2}\|\dot x(t)\|^2 = \mathcal{E}(t) -\frac{2\beta(t)}{\gamma+2}\langle x(t) - x^*,\dot x(t)\rangle,$$
we get:
\begin{eqnarray*}
\mathcal{E}'(t) &\leqslant & -\frac{2\gamma}{\gamma+2}\beta(t)\mathcal E(t)  + \frac{2\beta(t)}{\gamma+2}\left(\frac{\dot\beta(t)}{\beta(t)} +\frac{\gamma-2}{\gamma+2}\beta(t)\right)\langle x(t)-x^*,\dot x(t)\rangle+ \langle g(t),\dot x(t)+\frac{2\beta(t)}{2+\gamma}(x(t)-x^*)\rangle
\end{eqnarray*}
as expected.

\subsection{Proof of Lemma \ref{lemma HB}: differentiating the energy function of Theorem \ref{th:unperturbed:flat case}} \label{Appendix HB}

Let $\lambda$ be a non-negative real constant and $\xi(.)$ a real-valued function. Consider the energy:
\begin{eqnarray*}
\mathcal{E}(t) &=& t^{2}(F(x(t)) - F^*) + \frac{1}{2}\|\lambda(x(t) - x^*) + t\dot{x}(t)\|^2 + \frac{\xi(t)}{2} \|x(t) - x^*\|^2\\
&=& t(a(t)+b(t)+\xi c(t))
\end{eqnarray*}
where:
$$   a(t)  = t(F(x(t)) - F^*),~b(t) = \frac{1}{2t}\|\lambda(x(t) - x^*) + t\dot{x}(t)\|^2,~c(t)  = \frac{1}{2t}\|x(t) - x^*\|^2.$$

\begin{enumerate}
\item \textbf{Differentiating} $\mathcal{E}$
\begin{equation} \label{thm 4.1. energy step 2 }
\begin{aligned}
\mathcal{E}'(t)  = &2t(F(x(t)) - F^*) + t^{2} \langle \nabla F(x(t)), \dot{x}(t) \rangle + \langle \lambda (x(t) - x^*) + t \dot{x}(t), (1+\lambda)\dot{x}(t) + t \ddot{x}(t)\rangle \\
& + \frac{\dot\xi(t)}{2} \|x(t) - x^*\|^2 + \xi(t) \langle \dot{x}(t), x(t) - x^*\rangle 
\end{aligned}
\end{equation}
Since $x$ satisfies the ODE (\ref{ODE:HB}) with $g\equiv 0$, we have:
\begin{equation*}
\ddot{x}(t) = -\frac{\alpha}{t^\theta} \dot x(t) -\nabla F(x(t)),
\end{equation*}
hence:
\begin{equation*} 
\begin{aligned}
\mathcal{E}'(t)  = &~2a(t)- \lambda t\langle \nabla F(x(t)), x(t) - x^* \rangle+\left[\lambda(\lambda+1-\alpha t^{1-\theta})+\xi(t)\right]\langle  x(t) - x^*,\dot{x}(t)\rangle \\
& +t\left(\lambda+1-\alpha t^{1-\theta}\right)\|\dot x(t)\|^2 + \frac{\dot\xi(t)}{2} \|x(t) - x^*\|^2.
\end{aligned}
\end{equation*}
Since $F$ satisfies $\Hu(\gamma_1)$, we get:
\begin{equation}
\begin{aligned}
\mathcal{E}'(t)  = &~(2-\lambda \gamma_1)a(t)+\left[\lambda(\lambda+1-\alpha t^{1-\theta})+\xi(t)\right]\langle  x(t) - x^*,\dot{x}(t)\rangle \\
& +t\left(\lambda+1-\alpha t^{1-\theta}\right)\|\dot x(t)\|^2 + \frac{\dot\xi(t)}{2} \|x(t) - x^*\|^2.
\end{aligned}
\end{equation}

Noticing that 
\begin{equation*}
\|\lambda(x(t) - x^*) + t \dot{x}(t) \|^2 = t^2\|\dot{x}(t)\|^2 + 2t\lambda \langle \dot{x}(t), x(t) - x^* \rangle + \lambda^2 \|x(t) - x^*\|^2
\end{equation*}
i.e.
\begin{equation}
t \|\dot{x}(t)\|^2 = 2 b(t) - 2 \lambda^2 c(t) - 2\lambda\langle x(t) - x^*, \dot{x}(t) \rangle,
\end{equation}
we get:
\begin{equation}
\begin{aligned}
\mathcal{E}'(t)  =~& (2-\lambda\gamma_1)a(t) + 2(\lambda+1-\alpha t^{1-\theta}) b(t) -2\lambda^2(\lambda+1-\alpha t^{1-\theta})c(t)\\
&+ \left[\xi(t) - \lambda(1+\lambda-\alpha t^{1-\theta}) \right]\langle x(t)-x^*,\dot x(t) \rangle + t\dot\xi(t)c(t).
\end{aligned}
\end{equation}
Setting $\underline{\xi(t) = \lambda(\lambda + 1 - \alpha t^{1-\theta})}$, we obtain:
\begin{equation}
\begin{aligned}
\mathcal{E}'(t)  =~& (2-\lambda\gamma_1)a(t) + 2(\lambda+1-\alpha t^{1-\theta}) b(t) -2\lambda^2(\lambda+1-\alpha t^{1-\theta})c(t)\\
&+ \left[\xi(t) - \lambda(1+\lambda-\alpha t^{1-\theta}) \right]\langle x(t)-x^*,\dot x(t) \rangle -\lambda\alpha(1-\theta) t^{1-\theta}c(t).
\end{aligned}
\end{equation}

\item \textbf{Differentiating} $\mathcal{H}$.

Recall that $\mathcal{H}(t) = t^p \mathcal{E}(t)$. Hence $\mathcal{H}'(t) = t^{p-1}(p\mathcal{E}(t) + (1+t)\mathcal{E}'(t))$
Hence
\begin{eqnarray*}
\mathcal{H}'(t) & =& t^{p-1}(p\mathcal E(t) + t\mathcal E'(t))\\
&=&t^{p}(pa(t) + pb(t) + p\xi(t)c(t) + \mathcal{E}'(t)) \\
& \leqslant &t^p\left[(2-\lambda\gamma_1+p)a(t) + (2\lambda+2-2\alpha t^{1-\theta}+p)b(t)\right.\\
&&\left.+\lambda((\lambda+1)(p-2\lambda) - \alpha(p+1-\theta-2\lambda)t^{1-\theta})c(t) \right]
\end{eqnarray*}
as expected.

\end{enumerate}

\bibliographystyle{plain}
\bibliography{reference}

\end{document}